\documentclass[12pt]{amsart}
\usepackage{amssymb}
\usepackage{amsmath}
\newtheorem{thm}{Theorem}[section]
\newtheorem{cor}[thm]{Corollary}
\newtheorem{lem}[thm]{Lemma}
\newtheorem{notation}[thm]{Notations}
\newtheorem{prop}[thm]{Proposition}
\newtheorem{ex}[thm]{Example}
\newtheorem{exs}[thm]{Examples}
\newtheorem{defns}[thm]{Definitions}
\theoremstyle{definition}
\newtheorem{defn}[thm]{Definition}
\theoremstyle{remark}
\newtheorem{rem}[thm]{Remark}
\newtheorem{rems}[thm]{Remarks}

\numberwithin{equation}{section}

\newcommand{\de}{\delta}
\newcommand{\De}{\Delta}

\newcommand{\df}{{\it 2}-fir}
\newcommand{\dfs}{{\it 2}-firs}
\newcommand{\intdel}{\bigcap_{\delta \in \Delta}R\delta}
\newcommand{\intgam}{\bigcap_{\gamma \in \Gamma}R\gamma}
\newcommand{\find}{$F$-independence \;}
\begin{document}

\title
{Algebraic and $F$-Independent sets in 2-firs}%

\maketitle
%----------------------------------------------------------------

\begin{center}
{\bf Andr\'{e} Leroy$^\dagger$ and Adem Ozturk$^\ddagger$}%

\markboth{\bf A. LEROY AND A.OZTURK}{ \bf ALGEBRAIC AND
$F$-INDEPENDENT SETS }

\vspace{5mm}
 {$^\dagger$Universit\'{e} d'Artois,
 Facult\'{e} Jean Perrin\\
Rue Jean Souvraz, 62 307 Lens, France. \\
leroy@euler.univ-artois.fr} \vspace{3mm}

{$^\ddagger$Universit\'{e} de Mons-Hainaut,
        Institut de Math\'{e}matique \\
Avenue du Champ de Mars, B-7000 Mons, Belgique. \\
ozturk@umh.ac.be}
\end{center}

\thanks{ }%
%\keywords{}%

%\date{}%
%\dedicatory{}%
%\commby{}%
% ----------------------------------------------------------------
\begin{abstract}
 Let $R$ denote a \df.  The notions  of $F$-independence
and algebraic subsets of $R$ are defined.  The decomposition of an
algebraic subset into similarity classes gives a simple way of
translating the $F$-independence in terms of dimension of some
vector spaces.  In particular to each element $a \in R$ is
attached a certain algebraic set of atoms and the above
decomposition gives a lower bound of the length of the atomic
decompositions of $a$ in terms of dimensions of certain vector
spaces. A notion of rank is introduced and fully reducible
elements are studied in details.
\end{abstract}

\vspace{5mm}

\section{Introduction and preliminaries}

\vspace{7mm}

The main goal of this paper is to study factorizations in \dfs
\ via a careful use of classical notions such as similarity
and a systematic use of new notions such as algebraicity and
$F$-independence.  An attempt has been made to keep the paper
relatively self-contained and examples have been given all along
the paper to facilitate the reading.

Let us recall that a ring $R$ is a \df\ if any right ideal of $R$
generated by at most 2 elements is free of unique rank.  Of
course, a {\df} is a domain and it can be shown  (Cf. \cite{C1})
that this definition is equivalent to the following one : A domain
$R$ is a {\df} if and only if

$$\forall a,b \in R, \; aR \cap bR \ne 0 \Rightarrow \exists c,d
\in R \; : \; aR \cap bR = cR \; ; \; aR + bR = dR \,.$$

The lack of symmetry in this definition is only apparent and in
the paper we will freely use the fact that it is in fact
symmetric.  For the convenience of the reader we include a proof
of this fact.  We first state a useful lemma (Cf. \cite{C1} and
\cite{C2}).
\begin{lem}
\label{similarity}
 Let $R$ be a domain and
$a,a'$ be nonzero elements in $R$.  Then, the following are
equivalent :
  \begin{enumerate}
  \item[(i)] $R/aR \cong  R/a'R$.
  \item[(ii)] $\exists \; b \in R$ such that $aR+bR=R$ and $aR \cap
  bR=ba'R$.
  \item[(iii)] $\exists \; b' \in R$ such that $Ra' +Rb' =R$ and
  $Ra' \cap Rb' =Rab'$.
  \item[(iv)] $R/Ra \cong  R/Ra'$.
  \end{enumerate}

  If $b\in R$ is as in ii) above, there exists $b'\in R$
  satisfying the equalities in iii) and such that $ba'=ab'$.
  Moreover we then have $R/Rb \cong R/Rb'$.

\end{lem}
\begin{proof}
$i) \Longleftrightarrow ii)$ and $iii) \Longleftrightarrow iv)$
are easy and left to the reader.

$ii) \Longrightarrow iii)$  Since $ba' \in aR$ and $aR+bR=R$, one
can find $b',c',d'\in R$ such that $ba'=ab'$ and $ad'-bc'=1$. This
leads to $a(d'a-1)=bc'a \in aR \cap bR = ba'R =ab'R$.  Hence there
exists $c\in R$ such that $c'a=a'c$ and $d'a -1 = b'c$. Similarly,
we have $ad'b=b(c'b+1) \in aR \cap bR = ba'R$, hence there exists
$d \in R$ such that $c'b +1 = a'd$.  We now get $a'cb' =
c'ab'=c'ba' = (a'd-1)a'=a'(da'-1)$.  This gives $cb'=da'-1$ and
thus $Rb'+Ra'=R$.

Now, if $x=pa'=qb' \in Ra' \cap Rb'$ we get $q(d'a - 1) = qb'c =
pa'c = pc'a \in Ra$.  This shows that $q \in Ra$ and $x\in
Rab'=Rba'$.  We conclude $Ra' \cap Rb'=Rab'$.

$iii) \Longrightarrow ii)$ This is given by duality using the
opposite ring $R^{op}$.

The last statement can be obtained by finding the right equations
in the above proof and by using the following:
 $$
 \frac{R}{Rb^\prime}\cong \frac{Ra^\prime +
 Rb^\prime}{Rb^\prime}\cong \frac{Ra^\prime}{Ra^\prime \cap
 Rb^\prime}\cong \frac{Ra^\prime}{Rab^\prime}\cong
 \frac{Ra^\prime}{Rba^\prime}\cong\frac{R}{Rb} \,.
 $$
\end{proof}
\vspace{5mm}

Let us mention that the equivalence $(i) \longleftrightarrow (iv)$
is due to Fitting (Cf \cite{F}).  We now get the desired
left-right symmetry of the definition of a \df\:

\begin{cor}
\label{symmetry} Let $R$ be a domain.  The following are
equivalent :

  $i)$ $\forall a,b \in R, \, aR \cap bR \ne 0 \Rightarrow
\exists c,d \in R \; : \; aR \cap bR = cR \; ; \; aR + bR = dR \,
$.

  $ii)$ $\forall s,t \in R, \, Rs \cap Rt \ne 0 \Rightarrow
\exists u,v \in R \; : \; Rs \cap Rt = Ru \; ; \; Rs + Rt = Rv$.

\end{cor}
\begin{proof}
Of course, we will only prove that $i)$ implies $ii)$.  So let
$s,t \in R$ be such that $Rs \cap Rt \ne 0$.  We can find $a,b \in
R$ such that $0 \ne as = bt$ and $i)$ shows that there exist $c,d
\in R$ such that $aR \cap bR = cR$ and $aR + bR =dR$.  Writing
$c=ab'=ba' \;,\; a=dx$, and $b=dy$, we get $dxb' = ab' =ba'
=dya'$.  Since $R$ is a domain this gives $xb'=ya'$ and we easily
obtain that $xR \cap yR = xb'R=ya'R$ and $xR +yR =R$.  Lemma
\ref{similarity} shows that $Ra' + Rb' =R$ and $Ra'\cap Rb' =
Rxb'=Rya'$.  Now, since $as=bt\in aR\cap bR=cR$ there exists $v
\in R$ such that $as=bt=cv=ab'v=ba'v$ and so $s=b'v$ and $t=a'v$.
We thus get the desired conclusions : $Rs +Rt =Rv$ and $Rs \cap Rt
= Ru$ for $ u=xb'v$.
\end{proof}

Lemma \ref{similarity} and Corollary \ref{symmetry} will be used
several times. For more details on {\df} we refer to P.M. Cohn's
book "Free rings and their relations" (\cite{C1}).  We assume now
that $R$ is a {\df} and we will analyze injectivity and
surjectivity of some maps.  For $a,a'\in R \setminus \{0\}$, a
nonzero $R$-module homomorphism $\phi :R/Ra \longrightarrow R/Ra'$
is determined by an element $b' \in R \setminus Ra' $ such that
$\phi (x + Ra) = xb' + Ra'$ for any $x \in R$.  For the map $\phi$
to be well defined we must have $ab' \in Ra'$, and hence there
exists $b \in R$ such that

$$0 \ne ba' = ab'.$$

In particular this implies that there exists a right $R$-module
homomorphism : $\phi' : R/a'R \longrightarrow R/aR$ given by
$\phi'(y+a'R) = by + aR$ for any $y\in R$.  Notice that, since $R$
is a domain, $b'\notin Ra'$ implies that $b \notin aR$; this shows
that $\phi'$ is also nonzero.  The next lemmas will establish a
kind of duality between these two maps.

\vspace{4mm}
\begin{lem}
\label{injectivity}
 Let $R$ be a {\df} and $a,a' \in R\setminus\{0\}$.
With the above notations the following are
 equivalent:
\begin{enumerate}
\item[(i)] $\phi$ is injective.
\item[(ii)] $xb' \in Ra' \Longrightarrow x \in Ra$.
\item[(iii)] $Ra' \cap Rb' =Rab' =Rba'$.
\item[(iv)] $aR+bR=R$.
\item[(v)] $ \exists \;d' \in R $ such that $bd' -1 \in aR$.
\item[(vi)] $\phi'$ is surjective.
\end{enumerate}
\end{lem}
\begin{proof}

$(i)\Leftrightarrow (ii)$ This is obvious.

$(ii)\Leftrightarrow (iii)$  We always have $Rab^\prime =
Rba^\prime \subseteq Ra^\prime \cap Rb^\prime$.  On the other hand
if $d=xb^\prime\in Ra^\prime \cap Rb^\prime$ and (ii) holds, then
$ x\in Ra$ and $d\in Rab^\prime = Rba^\prime$.  Conversely if
(iii) holds and $xb^\prime \in Ra^\prime$, then $xb^\prime \in
Ra^\prime \cap Rb^\prime$.  Since $R$ is a domain we get $x\in
Ra$.

$(iii)\Leftrightarrow (iv)$  Assume $(iii)$ holds.  We have
$0\not=ab^\prime=ba^\prime \in aR \cap bR$, and, since $R$ is a
{\df}, we can write $aR\!+\!bR=dR$ for some $d \in R$.
 In particular, there exist $x,y \in R$ such that $a=dx$ and
 $b=dy$.  So $dxb' = ab' =ba'=dya'$ and we get $xb' \in Ra' \cap
 Rb' = Rab'$.  This gives $x \in Ra = Rdx$ and we conclude that
$d$ is a unit in $R$ and $aR+bR=R$.  Now, assume $(iv)$ holds.
Since $0\not=ab^\prime \in Ra^\prime \cap Rb^\prime$ we know that
there exists $x\in R$ such that $Ra^\prime \cap Rb^\prime =Rx$.
Let $r,s,t\in R$ be such that $x=sa^\prime=tb^\prime$ and
$ab^\prime =ba^\prime =rx$.  We then get $a=rt$ and $ b=rs$. Using
these equalities and $(iv)$ we get that $1=au +bv = rtu+rsv$.
Hence $r$ is a unit in $R$ which implies $Rab^\prime =Rrx=Rx=
Ra^\prime \cap Rb^\prime$, as desired.

The other equivalences are easy and left to the reader.
\end{proof}

As we have seen the notion of a {\df} is left-right symmetric,
hence using similar arguments as the ones used in the above proof
we also get the following:

\vspace{5mm}
\begin{lem}
\label{surjectivity}
 With the notations of the previous lemma, the following are equivalent:

\begin{enumerate}
\item[(i)] $\phi$ is surjective.
\item[(ii)] $ \exists \; d \in R $ such that $db^\prime -1 \in
Ra^\prime$.
\item[(iii)] $Ra^\prime + Rb^\prime=R$.
\item[(iv)] $aR \cap bR =ab'R =ba'R$.
\item[(v)] $bx \in aR \Longrightarrow x \in a'R$.
\item[(vi)] $\phi'$ is injective.
\end{enumerate}
\end{lem}
\begin{proof}
This is left to the reader.
\end{proof}

One of our aims in this paper is to analyze atomic factorizations
of elements of a {\df} $R$ using dimensions of some vector spaces
over division rings of the form $End_R(R/Rp)$ where $p$ is an atom
of $R$ \footnote{Let us recall that a nonzero element in a ring
$R$ is an atom if it is not a unit and cannot be written as a
product of two non units}.  If $R$ is left principal
$\frac{R}{Rp}$ is simple and Schur's lemma implies that
$End_R(R/Rp)$ is a division ring.  For an atom $p$ in a {\df} $R$
it is not true, in general, that $R/Rp$ is a simple module (Cf.
\ref{examples of 2-firs} d), below), nevertheless, as is well
known, the analogue of Schur's lemma is true.  We include a short
proof for completeness:

\vspace{5mm}
\begin{cor}
\label{Schur} Let $p$ be an atom in a {\df} $R$.  Then
$End_R(R/Rp)$ is a division ring.
\end{cor}
\begin{proof}
Let $\phi \in End_R(\frac{R}{Rp}) \setminus \{0\}$ and put $\phi
(1+Rp) = b' + Rp $. There exist $b,d\in R$ such that $0\not =
pb'=bp \in pR \cap bR$ and $pR + bR = dR$.  In particular, there
exists $d' \in R $ such that $p = dd'$.  Assume $d'$ is a unit
then $b \in dR = pR$.  We then get $pb'= bp\in pRp$, so that
$b'\in Rp$.  But this contradicts the fact that $\phi \not = 0$
and so $d'$ cannot be a unit.  Since $p=dd'$ is an atom we must
have that $d$ is a unit and $pR + bR=R$.  Lemma \ref{injectivity}
implies that $\phi$ is injective.  A similar argument shows that
$\phi':R/pR \longrightarrow R/pR$ defined by $\phi'(1+pR)=b+pR$ is
also injective and so Lemma \ref{surjectivity} implies that $\phi$
is surjective.
\end{proof}

\vspace{5mm}
\begin{rem}
A complete characterization of elements $a \in R$ which are such
that $End_R(R/Ra)$ is a division ring has been obtained by the
second author (\cite{O}).
\end{rem}
\vspace{9mm}

We close this section with some examples :

\vspace{3mm}

\begin{exs}
\label{examples of 2-firs} {\rm
\begin{enumerate}
\item[a)] A {\df} is a domain and, since in a commutative domain
$R$ two nonzero elements $a,b \in R$ are always such that $0\ne ab
\in aR \cap bR$, we see that a commutative {\df} is simply a
domain in which every finitely generated ideal is principal. These
rings are called B\'{e}zout domains in the literature.
\vspace{2mm}
\item[b)] In the same spirit, if $R$ is a right noetherian domain
then, as is well-known and easy to check, $0 \ne aR \cap bR $ for
any nonzero elements $a,b \in R$.  Hence a domain is a right
noetherian {\df} if and only if it is principal. Of course, a
{\df} which is an Ore ring is a Bezout domain but is not always
principal (Consider for example the commutative ring $R=\mathbb Z
+ x\mathbb Q[[x]]$ discussed in the next section, Cf example
\ref{non atomic 2-fir}). \vspace{2mm}
\item [c)] One good source of inspiration for our purpose is the
case of an Ore extension : $R=K[t;S,D]$ where $K$ is a division
ring, $S$ an endomorphism of $K$ and $D$ an $S$-derivation (Let us
recall that the elements of $R$ are polynomials
$\sum_{i=0}^{n}a_it^{i}$ with coefficients $a_i \in K$ written on
the left and the commutation rule is given by $ta = S(a)t + D(a)$)
.  Since $R$ is always a left principal ideal domain but $R$ is
right principal if and only if $S$ is onto, the ring $R$ is a
 {\df} which is not necessarily a right PID.  This ring and
factorization of its elements have been extensively studied in
\cite{LL},\cite{LL1} and \cite{LL2}.  These papers were starting
points for our reflections. \vspace{2mm}
\item[d)] Let $k$ be a field and $T=k(x)[t;S]$ be the Ore
extension where the $k$-endomorphism $S$ is defined by $S(x)=x^2$.
Consider $R=T^{op}$ the opposite ring of $T$.  $R$ is a {\df},
$t-x$ is an atom of $R$ but we claim that $\frac{R}{R(t-x)}$ is
not a simple $R$-module.  Equivalently we must show that
$\frac{T}{(t-x)T}$ is not a simple right $T$-module.  For any
$a\in k(x)$ there exist $a_0\, ,\; a_1 \in k(x)$, uniquely
determined, such that $a=S(a_0) + xS(a_1)$, and we have
$at=(t-x)a_0 + xta_1 +xa_0$.  From this it is easy to check that
$T = (t-x)T \bigoplus xtT \bigoplus k(x)$ (using induction on the
degree to get a decomposition of any polynomial and the fact that
$x \notin S(k(x))$ in order to prove that the sum is direct).  We
get that $\frac{T}{(t-x)T} = xtT \bigoplus k(x)$  (notice that the
right $T$-module structure of $k(x)$ is given by $a.t = xa_0$).
This module is obviously not simple and in fact it is not even
semisimple since it is finitely generated but not artinian.
\vspace{2mm}
\item[e)] If $k$ is a field, the free $k$-algebra $k<x,y>$ is
a {\df}.  We refer to P.M.Cohn's book for a proof of this fact
(Cf. \cite{C1}).

\end{enumerate}
}
\end{exs}

%---------------------------------------------------
%---------------------------------------------------
\vspace{10mm}

\section{Length and similarity}

\vspace{7mm}

Let us start with a definition which is crucial while dealing with
factorization in a non-commutative setting.

\vspace{5mm}
\begin{defn}
Two nonzero elements $a,a'$ in a domain $R$ are {\bf similar} if
$R/Ra \cong R/Ra'$.  We will then write $a \sim a'$.
\end{defn}
Lemma \ref{similarity} shows that this notion is left-right
symmetric and provides other characterizations of this definition.
In fact this notion defines an equivalence relation on the set
$R$.  The decomposition into similarity classes will play an
important role in our considerations.

\vspace{5mm}
\begin{exs}
{\rm
\begin{enumerate}
\label{examples of similarity}
\item[a)] Two elements $a,a' \in R\setminus\!\{0\}$ are associate
(resp. right associate or left associate) if there exist
invertible elements $u,v \in R$ such that $a'=uav$ (resp. take
resp. u=1 or v=1).  We leave to the reader to check that associate
elements are in fact similar.
\item[b)] In the case of an Ore extensions $R=K[t;S,D]$ where $K$
is a division ring and $D$ is an $S$-derivation, two elements $t
-a$ and $t-a'$ are similar if and only if there exists a nonzero
$c\in K$ such that $a'c = S(c)a + D(a)$.  This plays an important
role in the evaluation and in the factorization theories in Ore
extensions \cite{LL}, \cite{LL1} and \cite{LL2}.
\end{enumerate}
}
\end{exs}
In this section we want to investigate the relations between
similarity and length.  We will try to avoid assuming that our
{\df} is atomic.  Let us first give an example of a non atomic
{\df}.  This classical example (\cite{A}) will also be used later
in the paper.

\vspace{5mm}
\begin{ex}
\label{non atomic 2-fir}

{\rm We remarked in \ref{examples of 2-firs} that a commutative
ring is  a {\df} if and only if it is a B\'{e}zout domain. We will
show that the ring $R=\mathbb Z + x\mathbb Q[[x]]$ is a \df\ but
is not atomic. For a nonzero series
$f=\sum_{i=0}^{\infty}a_ix^{i}$ we define $o(f)=min\{i\in \mathbb
N \vert a_i \ne 0 \}$.  This series is invertible if and only if
$a_0 \in \{+1,-1 \}$.  Let us remark that any element $f$ of $R$
can be written in the form $a_mx^mu$ where $a_m\in \mathbb Q,
m=o(f), u \in U(R)$, the units of $R$ (obviously if $m=0$, $a_0
\in \mathbb Z)$.
 This implies that the nonzero principal ideals of $R$ are of the form
 $a_mx^mR$ where $a_m \in \mathbb Q$ if $m>0$ and $a_m \in \mathbb
 N$ if $m = 0$.  Let $A=a_mx^mR$ and $B=b_lx^lR$ be two nonzero principal
 ideals.  In order to show that $R$ is a \df \ we must prove that $A+B$
 is again principal.  Without loss of
 generality we may assume that $m \le l$.  We consider three cases:
\begin{enumerate}
\item[case 1]
If $m=l=0$ then $A=a_0R, \; B=b_0R$ and $A+B=dR$
 where $d$ is a greatest common divisor of $a_0$ and $b_0$ in
 $\mathbb Z$
 \item[case 2]
If $0\le m < l$ then $b_lx^lR =a_mx^mb_la_m^{-1}x^{\ell -m} R
\subset a_mx^mR$ hence $A + B = a_mx^mR$
 \item[case 3]
 If $0<l=m$.  Let us write $a_m=ce^{-1}, \; b_m = de^{-1}$
 with $c,d,e \in \mathbb Z$ and $e\neq 0$.  Define $f$ to be the least common multiple of $c$
 and $d$. Now we have
$a_mx^mR \cap b_mx^mR = (e^{-1}x^m)(cR\cap dR)= e^{-1}x^{m}fR$.  We conclude that $R$ is indeed
a \df .  Now, $2 \in R$ is obviously an atom and for all $n \in
\mathbb N$ we can write $x=2^n(x2^{-n})$. This shows that $x$ is
divisible by elements of arbitrary length. Hence $R$ cannot be
atomic.
\end{enumerate}
}
\end{ex}

\vspace{5mm}

The following lemma is the the starting point for establishing the
relations between atomic factorizations and lengths.  A proof of
it is given for instance in \cite{LL}.

\begin{lem}
\label{similar atoms} Let $p$ be an atom in a {\df} $R$.  If $q
\in R$ is similar to $p$, then $q$ is also an atom.
\end{lem}

Let us mention that the result is not true if $R$ is a domain
which is not a \df, once again an example can be found in
\cite{LL}.

\vspace{5mm}
\begin{lem}
\label{factofsim} Let $a$ and $a'$ be two similar elements in a
{\df}. If $a=bc$ then there exist $b',c'\in R$ such that $b \sim
b', c \sim c'$ and $a'=b'c'$.
\end{lem}
\begin{proof}
Let $\phi :R/a'R \longrightarrow R/aR$ be an isomorphism
determined by $\phi(1+a'R)=x+aR$.  Lemma \ref{similarity} shows
that there exists $x'\in R$ such that $xR + aR =R,\; Rx' + Ra' =
R$ and $ax'=xa'$ i.e. $bcx'=xa'$. Since $R$ is a {\df} we have
$Rcx' + Ra'= Rc'$, for some element $c'\in R$. In particular there
exist $r,b' \in R$ such that $cx'=rc',a'=b'c'$. We claim that $b
\sim b'$.  Let us put $\psi : R/b'R \longrightarrow R/bR : y + b'R
\mapsto xy + bR$. This map is well defined since $xb'c' = xa' =ax'
= bcx' = brc'$ which shows that $xb' = br \in bR \setminus \{0\}$.
Moreover we have $R = xR + aR \subseteq xR + bR$, hence $xR + bR =
R$.  Now Lemma \ref{injectivity} implies that $\psi$ is
surjective.  On the other hand, $0\ne xb'=br \in Rr \cap Rb' $
and, since $R$ is a {\df}, we have $Rr + Rb' = Rd$ for some $d \in
R$.  But then, $Rc' = Rcx' + Ra'= Rrc' + Rb'c' =  Rdc'$ which
shows that $d$ is a unit. Hence $Rr + Rb' = R$ and the lemma
\ref{surjectivity} implies that $\psi$ is injective and so, we
conclude that $b \sim b'$.  This proves the claim.  As $a'=b'c'$,
it remains to show that $c \sim c'$.

 Define a right $R$-morphism
$\Gamma :\frac{b'R}{a'R} \longrightarrow \frac{bR}{aR} :b'y + a'R
\mapsto xb'y +aR$.  This map is well defined since $xa' \in aR$.
We claim that $\Gamma$ is an isomorphism.  If $\Gamma(b'y +a'R)=0$
we get $xb'y \in xR \cap aR = xa'R$, where the last equality is
due to the injectivity of $\phi$. Since $R$ is a domain, this
gives $ b'y \in a'R$. This shows that $\Gamma$ is injective. Let
us now show that it is also surjective: it is enough to show that
$b + aR \in \Gamma(\frac{b'R}{a'R})$. Since $xR + aR =R$, we know
that there exist $u,v \in R$ such that $xu + av =1$. Consider
$\psi(ub +b'R)=xub + bR=(1-av)b + bR= b(1 - cvb)+ bR = bR$.  Since
$\psi$ is an isomorphism, we conclude that $ub \in b'R$ and
$\Gamma(ub + a'R)=xub + aR = b + aR$, as desired. This means that
$\Gamma $ is an isomorphism.  We conclude
$$\frac{R}{c'R} \cong\frac{b'R}{a'R} \stackrel{\Gamma}\cong
\frac{bR}{aR} \cong \frac{R}{cR} \, ,$$ as required.
\end{proof}

\vspace{5mm}
\begin{rem}
\label{product of similar} Let us notice that we can have $b \sim
b'$ but $bc \not \sim b'c$.  Let $\mathbb H$ denote the division
ring of real quaternions and let $R= \mathbb H [t]$.  Since $(j -
i)j = -i(j - i)$ we have that $t - j \sim t + i$ (Cf. example b)
in \ref{examples of similarity} ).  It can be shown that the only
non trivial monic right factor of the polynomial $f(t) = (t - j)(t
- i)$ is $t - i$ hence the module $R /Rf$ cannot be semisimple but
for $g(t) = (t+i)(t-i)=t^2 +1$ the left $R$-module $R/Rg \cong
R/R(t-i) \bigoplus R/R(t-j)$ is semisimple.  This shows that
$f(t)$ is not similar to $g(t)$.
\end{rem}

\vspace{5mm}
\begin{defn}
Let $a$ be a nonzero element in a \df\ $R$.  An element $a$ in $R$
is of finite length $n$ if it can be factorized into a product of
$n$ atoms but does not admit a factorization into a product of
less than $n$ atoms . If an element in $R$ cannot be factorized
into a finite product of atoms we will say that it is of infinite
length.  The invertible elements of $R$ are of length $0$.
\end{defn}
We denote the subset of elements that can be written as products
of atoms by $\mathcal F = \{x \in R \vert \ell (x) < \infty \}$.
Notice in particular that units of $R$ are included in $\mathcal
F$, but $0 \notin \mathcal F$.

\vspace{5mm}
\begin{thm}
\label{length of similar elements}

In a {\df}, similar elements have the same length.
\end{thm}

\begin{proof}
Let $a,a'$ be nonzero similar elements in $R$.  We proceed by
induction on the length of $a$.  The claim is obvious if
$\ell(a)=0$.  The above lemma \ref{similar atoms} shows that the
theorem is true for $\ell(a)=1$.
   Now, let $n\ge 2$ and assume that the theorem is true
for elements of length $\leq n\!-\!1$ and let $a\in R$ be such
that $\ell(a)=n$.  Obviously we must have $\ell(a')\geq n$.  Let
us write $a=bq$ where $\ell(b)=n\!-\!1$.  So $q$ is an atom. Since
$a$ is similar to $a'$, lemma \ref{factofsim} shows that there
exist
 $b', \, q'$ such that $a' = b'q'$ where $q' \sim q$ and $b' \sim
 b$.   The induction hypothesis implies that $\ell(q')=1$ and
 $\ell(b')=n\!-\!1$.  We thus conclude that $\ell(a')=n=\ell(a)$.
 This proves the theorem.
\end{proof}

For two elements $a,b$ in a \df\ $R$ such that $Ra \cap Rb \ne 0$,
we will denote by $a^b$ an element in $R$ such that $Ra \cap Rb =
Ra^bb$.  Notice first that $a^b$ is defined up to a left multiple
by a unit.  We will also use $Ra\cap Rb=R[a,b]_\ell$ with
$[a,b]_\ell$ in $R$.  In the next lemma we briefly study some
properties of $a^b$.  Recall that for $(a,b) \in R^2 , \,
(Ra)b^{-1}= \{x \in R \vert xb \in Ra \}$.

\vspace{5mm}
\begin{lem}
\label{a^b}
\begin{enumerate}
\item [(a)]Let $a$ be an atom in $R$ and $b \in R\setminus Ra$ such
that $Ra\cap Rb\neq 0$.  Write $Ra\cap Rb=Ra^{b}b$. Then $a^b$ is an
atom similar to $a$.
Conversely : if $a$ is an atom and $a'\in R$ is such that $a'\sim
a $ then there exists $b\in R\setminus Ra$ such that
$Ra'=Ra^b$.  Moreover we have $Ra^b = (Ra)b^{-1}$.

\item [(b)]If $a,b,c$ are elements of $R$ such that $Ra \cap Rcb
\ne 0$ then we have $Ra^{cb} = R(a^b)^c$.
\item[(c)]If $Ra \cap Rb \cap Rc \ne 0$ then we have $R[a,b]_\ell^c=
R[a^c,b^c]_\ell$.
\end{enumerate}
\end{lem}

\begin{proof}(a) As $b\not\in Ra$, $a^b$ is not a unit. Let us remark
that $Ra\cap Rb\neq 0$ implies that $Ra+Rb=Rd$ for some $d\in R$.
As $a$ is an atom and $b\not\in Ra$, $d$ is  a unit in $R$. So, by
lemma \ref{similarity} \ $a^b \sim a$.  The converse and the
additional statement are left to the reader.

  (b) We have $Ra \cap Rb = Ra^bb$ and $Ra^b \cap Rc =
R(a^b)^cc$ (both of these intersections are nonzero thanks to our
assumption in (b)). Multiplying the last relation by $b$ on the
right we get $R(a^b)^ccb = Ra^bb \cap Rcb$. Using the first
relation this leads to $R(a^b)^ccb = Ra \cap Rb \cap Rcb =Ra \cap
Rcb = Ra^{cb}cb$ and this gives the desired relation.

  (c) $R[a,b]_\ell^cc = R[a,b]_\ell \cap Rc = Ra \cap Rb
\cap Rc = Ra \cap Rc \cap Rb \cap Rc = Ra^cc \cap Rb^cc =
R[a^c,b^c]_\ell c$.  This leads to the desired equality.

\end{proof}

\vspace{5mm}
\begin{lem}
\label{product formula} Let $a,b$ be nonzero elements of $R$ and
$p$ an atom in $R$ such that $ab \in Rp$ but $b \notin Rp$.  Then
$a \in Rp^b$.
\end{lem}
\begin{proof}
Since $0 \ne ab \in Rp \cap Rb$, we know that $Rp \cap Rb =
Rp^bb$.  In particular there exists $c \in R $ such that $ab =
cp^bb$ and we get $a = cp^b$.
\end{proof}

Since atomic factorizations of two elements $b$ and $c$ yield an
atomic factorization of their product $bc$, we have $\ell(bc) \leq
\ell(b)+\ell(c)$. The reverse inequality is not completely clear
and we offer a short proof of it in the next lemma.

\vspace{5mm}
\begin{lem}
\label{factor of atomizable} Let $R$ be \df.  Then
 $$ \ell(bc) = \ell(b) + \ell(c) \;\;\ \forall \;\; b,c \in R $$
\end{lem}
\begin{proof}
The case when $a:=bc$ has infinite length is clear and we may
assume that $\ell(a)=n<\infty$.  We must show that $n = \ell(b) +
\ell (c) $.  We proceed by induction on $n$.  The claim is obvious
for $n=0$.  If $n=1$, $a$ is an atom and the result is clear.
Assume $n>1$, obviously we may assume that neither $b$ nor $c$ are
invertible.  Write $a=p_1p_2 \cdots p_n = bc $. If $c \in Rp_n$
then there exists $c'\in R$ such that $c = c'p_n$ and we get
$bc'=p_1p_2 \cdots p_{n-1} $ and the induction hypothesis allows
us to conclude easily.  We may thus assume that $c \notin Rp_n $.
We have

$$ a=p_1p_2 \cdots p_n = bc \in Rc \cap Rp_n=Rc'p_n=Rp_n'c$$

 where $p_1,p_2, \dots ,p_n$ are given atoms and
$c',p_n'$ are elements in $R$.   Notice also that by Theorem \ref
{length of similar elements} and Lemma \ref{a^b} we have $\ell
(p_n') = \ell (p_n) = 1$ and $\ell (c') = \ell (c) \ge 1 $ .  The
above displayed equality  shows that there exist $r \in R$ and
$\alpha$ a unit in $R$ such that $bc = a =rp_n'c$ and $p_n'c =
\alpha c'p_n$ hence $b=rp_n'$ and $p_1p_2 \cdots p_{n-1} = r\alpha
c'$.  The induction hypothesis then shows that  we have $n-1 =
\ell (r) + \ell (c')$.  Since $b=rp_n'$ and $p_n'$ is an atom we
get $\ell(b) \le \ell (r) + 1 = n-1 - \ell (c') +1 \le n-1 $.
Hence we can again apply our induction hypothesis and we obtain
$\ell (b) = \ell (rp_n') = \ell (r) + \ell (p_n') = \ell (r) + 1 $
This gives $n = \ell (r) + \ell (c') + 1 = \ell (b) + \ell (c)$,
as desired

\end{proof}

The next theorem is part of folklore.

\vspace{5mm}
\begin{thm}
\label{length of gcd and llcm} Let $R$ be a \df\ and let $a,b \in
R\setminus\{0\}$ such that $Ra\cap Rb \neq 0$.  Write $Ra\cap
Rb=R[a,b]_\ell$ and $Ra+Rb=R(a,b)_r$. Then $$
\ell([a,b]_\ell)+\ell((a,b)_r)=\ell(a)+\ell(b)$$
\end{thm}

\begin{proof} If $\ell(a)$ or $\ell(b)$ are infinite then the equality is obvious.
The Noether's Isomorphism Theorem gives an isomorphism of
$R$-modules $(Ra+Rb)/Ra \cong Rb/Ra\cap Rb$. Let us write
$a=a''(a,b)_r$ for some $a'' \in R$ and $[a,b]_\ell=a'b$ with $a'
\in R$.  We get
 $$R/Ra'' \cong R/Ra' \; \, \mbox{i.e.}  \,\; a'' \sim a'$$

   So, by theorem \ref{length of similar elements}, $\ell(a'')=\ell(a')$. Now we
have
$\ell([a,b]_\ell)=\ell(a'b)=\ell(a')+\ell(b)=\ell(a'')+\ell(b)=
\ell(a)-\ell((a,b)_r)+\ell(b)$. This completes the proof.
\end{proof}

\vspace{5mm}
\begin{cor}
\label{length of a^b} Let $a,b$ be nonzero and non unit in $R$
such that $Ra\cap Rb \neq 0$.  Write $Ra\cap Rb=Ra^{b}b$.  Then
$\ell(a^{b}) \leq\ell(a)$.
\end{cor}
\begin{proof} The previous theorem gives $\ell(a^{b})+\ell(b)=\ell(a^{b}b) \leq
\ell(a)+\ell(b)$.
\end{proof}

\vspace{5mm} Let us now offer a few easy but important facts about
the subset $\mathcal F=\{x \in R \vert \ell (x) < \infty \}$ which
was introduced in the paragraph before \ref{length of similar
elements}.
\begin{prop}
\label{atomic elements} Let $a,b$ be elements in $\mathcal F $.
Then :
\begin {enumerate}
\item[a)] $ab \in \mathcal F$.
\item[b)] If $a \sim a'$, then $a' \in \mathcal F$.
\item[c)] If $Ra + Rb = Rc$, then $c \in \mathcal F$.
\item[d)] If $Ra \bigcap Rb  = Rd \ne 0$, then $d \in \mathcal F$.
\item[e)] If $\,\Gamma $ is a finite subset of $\mathcal F$ such that
$\bigcap_{\{\gamma \in \Gamma\}}R\gamma \ne 0$, then there exists
$h \in \mathcal F$ such that $\bigcap_{\{\gamma \in
\Gamma\}}R\gamma = Rh$.
\end{enumerate}
\end{prop}
\begin{proof}
  a) This is clear from \ref{factor of atomizable}.

  b) This comes from \ref{length of similar elements}.

  c) Since $a \in Rc$ this also follows from \ref{factor
of atomizable}.

  d) This follows from \ref{length of gcd and llcm}.

  e) This is obtained by repeated applications of the point e)
  above.

\end{proof}

We end this section with the following remark:

\begin{rem}
If $\Delta \subseteq R$ is such that $\intdel \bigcap \mathcal F
\ne \emptyset$ then $\Delta \subseteq \mathcal F$ \
\end{rem}
\begin{proof}
Let $x\in (\bigcap_{\delta \in \Delta}R \delta)\bigcap \mathcal
F$.  Then $\ell(x)<\infty$ and, since $x\in R\delta$ for any
$\delta \in \Delta$, we have that $\delta \in \mathcal F$ for all
$\delta \in \Delta$.
\end{proof}

\vspace{10mm}
%-----------------------------------------------------------------
%------------------------------------------------------------------

\section{$F$-independence on 2-fir}

\vspace{7mm}

We now  introduce some central definitions. In this section $R$
will denote a \df.  Due partly to the last remark and to make life
easier, we will only consider, in the next definition, subsets of
$\mathcal F$.

\vspace{5mm}
\begin{defn}
\label{F-algebraic} Recall that $U(R)$ stands for the set of
invertible elements of $R$.  A subset $\Delta$ of $\mathcal F
\setminus U(R)$ is said to be {\bf $F$-algebraic} if
$\bigcap_{\delta \in \Delta}R\delta \ne 0$.
\end{defn}

\vspace{5mm}

 Let us first remark that, according to the above definition,
 the empty set is algebraic (we follow the convention
 that the intersection of an empty family of subsets of $R$ is R itself).
 We now make some more remarks and introduce some notations
in the following:
\begin{notation}
{\rm If $\De$ is an $F$-algebraic subset of $R$ we will denote by
$\De_\ell$ an element (when there exists one) such that $\intdel =
R\De_\ell$.  For convenience and in accordance with future
notations,  we will put $\emptyset_\ell = 1$. If it exists
$\De_\ell$ is not unique but all such elements are left associates
and have the same length. We will sometimes use the word
"algebraic" meaning in fact "$F$-algebraic}".
\end{notation}

\begin{rems}
\label{algebraic rem}

 a) We have excluded the invertible elements from algebraic sets.
 The first reason is that algebraic sets are in fact a tool
 for the study of factorization the second reason is that if one
 admits invertible elements in algebraic sets this creates
 technical problems and more complicated statements.

b) Notice that if $\De \subseteq \mathcal F$ is algebraic and
finite, then, since $R$ is a \df,  $\De_\ell$ always exists and is
nonzero.  Moreover, Proposition \ref{atomic elements} e) shows
that in this case $\De_\ell \in \mathcal F$.  In particular, for
any finite subset $\Gamma$ contained in an algebraic set $\De$
there exists $\Gamma_\ell \in \mathcal F$ such that $\intgam =
R\Gamma_\ell$.

c) Remark also that if $\De$ is algebraic and $\De'$ is a subset
of $R$ consisting of right non invertible divisors of elements of
$\De$ then $\De'$ is also algebraic.

d) Let us mention that although $\De$ is a subset of $\mathcal F$,
the element $\De_\ell$, when it exists, might be of infinite
length (cf example \ref{algebraic example}, e).

e) In \cite{LL3} $(S,D)$-algebraic sets are defined in the context
of an Ore extension $R:=K[t;S,D]$ over a division ring $K$.  The
relation between this notion and the notion of $F$-algebraic sets
introduced above is as follows : a subset $\De \subseteq K$ is
$(S,D)$-algebraic if and only if the set $\{\, t - \de \,\vert \;
\de \in \De \, \} \subset R$ is $F$-algebraic in the sense defined
in \ref{F-algebraic}.

f) An $F$-algebraic subset of $\mathcal F$ should be called left
$F$-algebraic.  A similar definition for right $F$-algebraic sets
can be given.  Singleton sets of $\mathcal F$ not contained in
$U(R)$ are, of course, left and right $F$-algebraic but there are
sets with only 2 elements that are left $F$-algebraic but not
right $F$-algebraic : This is the case of $\{t,at\}\subset R =
k[t;S]$ where $k$ is a field, $S$ is an endomorphism of $k$ which
is not an automorphism and $a\in k \setminus S(k)$.  In this paper
$F$-algebraic will always refer to the left notion defined above.

\end{rems}
In the following $\mathcal A$ will stand for the set of atoms of
$R$. \vspace{4mm}
\begin{prop}
\label{equivalence for independence} Let $\Delta \subseteq
\mathcal F$ and $d\in \mathcal F$ be such that $\Delta \cup \{d \}
$ is $F$-algebraic.  Then the following are equivalent:
\begin{enumerate}
\item[(i)] There exist a finite subset $\Gamma$ of $\De$ and
$p \in \mathcal A$  such that $d \in Rp$ and $Rp + \intgam \ne R$.
\item[(ii)] There exists a finite subset $\Gamma$ of $ \Delta$ such that
$Rd + \bigcap_{\gamma \in \Gamma}R \gamma \ne R$.
\item[(iii)] There exist a finite subset $\Gamma$ of $ \Delta$ and
an element $\Gamma_\ell \in \mathcal F$ such that $\bigcap_{\gamma
\in \Gamma}R\gamma = R\Gamma_\ell$ and $R\Gamma_\ell \bigcap Rd =
Rm \not = 0$ with $\ell (m) < \ell (\Gamma_\ell) + \ell (d)$.
\end{enumerate}
\end{prop}
\begin{proof}
$(i)\Longrightarrow(ii)$ If $d ,p $ and $\Gamma$ are as in (i)
then $Rd + \intgam \subseteq Rp + \intgam \ne R$.
%As noticed in
%Remark \ref{algebraic rem} (b) above, there exists $\Gamma_\ell
%\in R$ such that $\intgam = R\Gamma_\ell$ and since $\Delta \cup
%\{d \} $ is {\it F-} algebraic, we have that $Rd \cap R\Gamma_\ell
%\ne 0 $. Hence (i) and Proposition \ref{atomic elements} show that
%there exists $f\in \mathcal F$ such that $Rd + R\Gamma_\ell = Rf
%\ne R$. If $p \in \mathcal A$ is a right factor of $f$ then $p$
%also divides both $d$ and $\Gamma_\ell$ on the right and we have
%$Rp + \intgam = Rp + R\Gamma_\ell = Rp \ne R$, as desired. For the
%converse implication

$(ii)\Longrightarrow (iii)$ Since $\Gamma$ is a finite subset of
$\mathcal F$, remark \ref{algebraic rem} b) shows that there
exists $\Gamma_\ell \in \mathcal F$ such that $\bigcap_{\gamma \in
\Gamma} R\gamma = R\Gamma_\ell $ and, since $\{d\} \cup \Gamma
\subseteq \{d\} \cup \Delta$ is algebraic, we have $0 \ne Rd \cap
R\Gamma_\ell =Rm$, for some $m$ in $R$. Since $R$ is a \df,
$Rd+R\Gamma_\ell=Ra$ for some $a$ in $R$. Now, \ref{atomic
elements} c) and d) give us that $a,m\in \mathcal F$.  From
Theorem \ref{length of gcd and llcm} we get $\ell (m) + \ell (a) =
\ell (d) + \ell (\Gamma_\ell) $. Since, by (ii), $a$ is not a unit
we have $ \ell (m) < \ell (d) + \ell (\Gamma_\ell)$ as required.

$(iii)\Longrightarrow (i)$ Since $R$ is a \df \ and  $R\Gamma_\ell
\bigcap Rd = Rm \not = 0$ there exists $a\in R$ such that
$R\Gamma_\ell + Rd = Ra$ and we have $\ell(d) +\ell(\Gamma_\ell)=
\ell(a)+\ell(m) < \ell(a)+\ell(d)+\ell(\Gamma_\ell)$.  Hence
$\ell(a)\not = 0$. If $p\in \mathcal A$ divides $a$ on the right
we have $d\in Rp$ and $Rp + \cap_{\gamma \in \Gamma}R\gamma = Rp +
R\Gamma_\ell \subseteq Rp + Ra = Rp \not = R$, as required.
\end{proof}
\vspace{3mm}

In view of the above proposition the following definitions appear
naturally:

\vspace{5mm}
\begin{defn}
\label{F-dependence} Let $\Delta$ be an algebraic set.
\begin{enumerate}
\item[(a)] An element $d \in \mathcal F$ is said to be {\bf $F$-dependent}
on $\Delta$ if $\De \cup \{d\}$ is algebraic and one of the
conditions of the above proposition is satisfied.
\item[(b)] $\Delta$ is {\bf $F$-independent}
if and only if for any $\delta \in \Delta$, $\delta$ is not
$F$-dependent over $ \Delta \setminus \{\delta \}$.
\item[(c)] An $F$-independent subset $B\subseteq\Delta$ is an
{\bf  $F$-basis} if any element of $\Delta$ is $F$-dependent on
$B$.
\end{enumerate}
\end{defn}

\medskip

It is clear that a subset of an $F$-algebraic set is also
$F$-algebraic and a subset of an $F$-independent set is also
$F$-independent.  The above proposition \ref{equivalence for
independence}(i) shows that it is possible to express
$F$-dependence by means of atoms.  We will have more precise
information in the next section.  For the moment let us notice the
following special case.
\begin{prop}
\label{dependence of an atom} Let $p$ be an atom in $R$.  Then $p$
is $F$-dependent on an $F$-algebraic set $\Delta \subset \mathcal F$ if and
only if there exists a finite subset $\Gamma \subseteq \Delta$
such that $\intgam \subseteq Rp$.
\end{prop}

The next proposition connects \find and length.

\vspace{4mm}
\begin{prop}
Let $R$ be a \df\ and $\Delta \subset \mathcal F$ be an $
F$-algebraic set.  Then $\Delta$ is an $F$-independent set if and
only if for any finite subset $\Gamma := \{a_1,a_2, \dots ,a_n \}
\subseteq \Delta$ we have $\ell (\Gamma_\ell ) =
\sum_{i=1}^{n}\ell (a_i)$ where $\Gamma_\ell \in \mathcal F$ is
such that $\cap_{i=1}^{n}Ra_i = R\Gamma_\ell$.
\end{prop}
\begin{proof} If some $a_j \in \{a_1,\ldots,a_n \}$ is
$F$-dependent on $\Gamma_j \!:=\!\{a_1,\ldots,a_n
\}\setminus\{a_j\}$, then $ \ell(\Gamma_\ell) <
\ell(a_j)+\ell((\Gamma_j)_\ell) \le \sum_{i=1}^{n}\ell(a_i)$,
where the last equality comes from Theorem \ref{length of gcd and
llcm}. Conversely, suppose $\{a_1,\ldots,a_n \}$ is
$F$-independent and let $a_j \in \{a_1,\ldots,a_n \}$.  By
induction, we may assume that $\ell((\Gamma_j)_\ell)=\sum_{i\neq
j}^{n}\ell(a_i)$. Since $a_j$ is not $F$-dependent on $\Gamma_j\!
=\!\{a_1,\ldots,a_n \}\setminus\{a_j\}$,
$\ell(\Gamma_\ell)=\ell(a_j)+\ell((\Gamma_j)_\ell)=\sum_{i=1}^{n}\ell
(a_i)$. This finishes the proof.
\end{proof}

\vspace{4mm}

\begin{cor}
\label{finite algebraic} Let $\Delta \subseteq \mathcal F $ be an
$F$-independent algebraic set in a \df \ $R$.  Then $\vert \De
\vert < \infty$ if and only if there exists $\De_\ell \in \mathcal
F$ such that $\intdel = R\De_\ell$.  Moreover in this case we have
$\vert\De \vert \le \ell(\De_\ell)$ and the equality occurs if and
only if $\De \subseteq \mathcal A$.

\end{cor}

  The above properties are quite nice but we will soon see
that the definitions of  $F$-dependence and  $F$-independence have
also some drawbacks.

There are some relations between $F$-bases and maximal
$F$-independent sets. To understand more precisely the
relationship, we first prove the following intermediate fact.

\vspace{4mm}
\begin{prop}
\label{exchange} Let $\Delta\cup\{a\}\cup\{b\}\subseteq \mathcal
F$ be an $F$-algebraic set in a \df \  $R$.  Then, if $b$ is
$F$-dependent on $\Delta\cup\{a\}$, but is not $F$-dependent on
$\Delta$, then $a$ is $F$-dependent on $\Delta\cup\{b\}$.
\end{prop}
\begin{proof} Since $b$ is $F$-dependent on
$\Delta \cup \{a\}$ there exists a finite subset $\Gamma$ of
$\Delta$ such that for $g \;,\; h \;,\; m \in \mathcal F$ defined
by $\bigcap_{\gamma\in\Gamma}R\gamma=Rg \; ,\; Ra\cap Rg = Rh$ and
$Rb\cap Rh=Rm$ we have $\ell (m)<\ell(h) + \ell(b) $. On the other
hand the fact that $b$ is not $F$-dependent on $\Delta$ implies
that $\ell(c) = \ell(b) + \ell(g)$ where $c$ is such that $Rb \cap
Rg = Rc$.  We thus have $\ell(a) + \ell(c) = \ell(a) + \ell(b) +
\ell(g) \ge \ell(b) + \ell(h) > \ell(m)$. Since we also have that
$Rm=Rb \cap Ra \cap Rg = Ra \cap Rc$ we conclude that $a$ is
$F$-dependent on $\Delta\cup\{b\}$, as required.
\end{proof}

\vspace{5mm}
\begin{prop}
\label{basis} Let $\Delta$ be an $F$-algebraic set in a \df \ $R$.
Then $B\subset\Delta$ is an $F$-basis of $\Delta$ if and only if $B$
is a maximal $F$-independent subset of $\Delta$. In particular,
any algebraic set has a basis.

\end{prop}

\begin{proof} The only if part is clear.  Assume $B$ is a maximal
$F$-independent subset of $\Delta$ and let $b\in \Delta\setminus
B$.  By assumption $B\cup\{b\}$ is not an $F$-independent set.
Hence some element $c\in B\cup\{b\}$ is $F$-dependent on the
others. If $c=b$, $b$ is $F$-dependent on $B$ as desired. Assume
$c\in B$.  Then $c$ is not $F$-dependent on $B\setminus\{c\}$ but
$F$-dependent on $(B\setminus\{c\})\cup \{b\}$.  By the last
proposition, $b$ is $F$-dependent on $(B\setminus\{c\})\cup
\{c\}=B$ as desired.  The last statement follows by using Zorn's
lemma.
\end{proof}

\vspace{5mm}
\begin{exs}
\label{algebraic example}
\begin{enumerate}
{\rm
\item[a)] Let $R=K[t;S,D]$ be an Ore extension over a division
ring $K$ where $S\in End(R)$ and $D$ is an $S$-derivation.  Let
$\{a_1, \ldots ,a_n \} $ be a subset of $K$ and consider $\Delta
:= \{t - a_1,\dots ,t - a_n \}\subset R$.  Then $\De$ is algebraic
since the $(t-a_i)$'s have a nonzero least common left multiple.
In fact, in this case, $R$ is a left principal ideal domain and
any finite subset of $R$ is (left) algebraic.  These situations
have been studied extensively in \cite{LL},\cite{LL1} and
\cite{LL2}. In these papers a basis for an algebraic set $\De$ was
called a P-basis.

\item[b)] Of course, a basis of an algebraic subset $\Delta$ of $R$ might
well be infinite.  When a basis is finite there exists an element
$h \in \mathcal F$ such that $ \bigcap_{b \in B} Rb=Rh$.  But in
general there may be no element $g\in R$ such that $\intdel =Rg$

\item[c)] Let us consider $R=K[t]$ where
$K$ is a field.  Let $a,b$ be nonzero element of $K$ and
$\Delta=\{(t-a)(t-b),(t-b)\}$. So $\Delta_\ell=(t-a)(t-b)$.
Moreover $B=\{t-b\}$ and $B'=\{(t-a)(t-b)\}$ are $F$-bases of
$\Delta$. But we have $R\Delta_\ell=RB'_\ell\subsetneq RB_\ell$.
This shows that even when an algebraic set $\De$ is finite the
least left common multiple of the element of a basis and the least
left common multiple of the elements of $\De$ may be different. It
will be shown later (Cf. Proposition \ref{B_l et D_l}) that such a
situation is impossible in the case when all elements of $\Delta $
are atoms.

\item[d)] Let $p_1,p_2$ be different atoms in $R$ such that $Rp_1 \bigcap
Rp_2 = Rm \ne 0$.  Then the set $\De = \{p_1,p_2,m\}$ is an
algebraic set.  Notice that $\{m\}$ and $\{p_1,p_2\}$ are bases
for $\De$ with different cardinals.

\item[e)] Let us consider the \df \ $R= \mathbb Z +
x\mathbb Q[[x]]$ (Cf. the example \ref{non atomic 2-fir}) and let
$\De =\{p \in \mathbb Z \vert {\rm p \; is \; prime \; and} \; p>0
\}$.  Notice that the elements of $\Delta$ are atoms in $R$.
Since $x \in \intdel $ we see that $\De$ is an algebraic subset of
$R$. In fact $\De$ is a basis of itself. Notice also that
$\intdel=Rx$. Since $x \notin \mathcal F$, this gives the example
promised in remark d) of \ref{algebraic rem}. }

\end{enumerate}
\end{exs}
\vspace{3mm}

The notions of $F$-dependence and $F$-independence are strongly
related to the notion of abstract dependence. Let us recall this
definition (Cf. \cite{J}).
\newline For a non vacuous set $X$ and a
relation $\Gamma$ from $X$ to the power set ${\mathcal{P}}(X)$, we
write $x\prec S$ if $(x,S)\in \Gamma$. We call $\Gamma$ a {\it
dependence relation} in $X$ if the following conditions are
satisfied :
\begin{enumerate}
     \renewcommand{\labelenumi}{(\roman{enumi})}
    \item if $x\in S$, $x\prec S$.
    \item if $x\prec S$, then $x\prec F$ for some finite subset
          $F\subset S$.
    \item if $x\prec S$ and every $y\in S$ satisfies $y\prec T$, then
          $x\prec T$.
    \item if $x\prec S$ but $x \nprec S\setminus\{y\}$ then $y
          \prec (S\setminus\{y\})\cup\{x\}$.
\end{enumerate}

   In our case $X=\mathcal F\setminus U(R)$, $S$ is an $F$-algebraic
set of $R$ and the relation ``$\prec$ '' is the $F$-{\it
dependence relation}. Obviously $(i)$ and $(ii)$ are satisfied.
The assertion $(iv)$ is given by \ref{exchange}.  But $(iii)$ is
false in general as the following example shows.

\vspace{5mm}
\begin{ex}
\label{no transitivity} {\rm Let $a,b,c,d$ be atoms in a \df \ $R$
such that $a$ is not similar to $d$ but $ba=cd$.   We then have
that $a$ is $F$-dependent on $\{ba\}$ and $ba=cd$ is $F$-
dependent on $\{d\}$ but $a$ is not $F$-dependent on $\{d\}$.}

\end{ex}

The problem of non transitivity disappears if we restrict
ourselves to algebraic sets of atoms.  Let us recall that
$\mathcal A $ denotes the set of atoms in $R$.

\vspace{5mm}
\begin{prop}
\label{transitivity for atoms} Let $\De,\;\De'\subseteq \mathcal
A$ be algebraic sets of atoms.  Assume $p\in \mathcal A$ is
$F$-dependent on $\De$ and each element of $\De$ is $F$-dependent
on $\De'$.  Then $p$ is $F$-dependent on $\De'$.
\end{prop}
\begin{proof}
By hypothesis there exists a finite subset $\Gamma$ of $\De$ such
that $\intgam \subseteq Rp$ (Proposition \ref{dependence of an
atom}).  Now each $\gamma \in \Gamma$ is $F$-dependent on $\De'$
and since $\Gamma$ is finite we can find a finite subset $\Gamma'$
of $\De'$ such that $\bigcap_{\gamma' \in \Gamma'}R\gamma'
\subseteq R\gamma$ for all $\gamma\in \Gamma$.  This means that
$\bigcap_{\gamma' \in \Gamma'}R\gamma' \subseteq \intgam \subseteq
Rp$.  This shows that $p$ is $F$-dependent on $\De'$.
\end{proof}

So if we restrict to algebraic sets of atoms the notion of
$F$-dependence defines an abstract dependence relation.  In this
case the general theory shows that a subset $B$ of an algebraic
set $\De$ is a basis if and only if it is minimal such that all
elements of $\De$ are $F$-dependent on $B$.

The restriction to subsets of $\mathcal A$ is not as bad as it
could seem on the first sight.  We have already seen that atoms
appear naturally while dealing with $F$-independence (see
\ref{equivalence for independence} (i)).  In Proposition
\ref{general dependence} we will show more precisely how the
notion of $F$-dependence on elements of $\mathcal F$ is controlled
by the $F$-dependence on $\mathcal A$.

\vspace{10mm}
%-------------------------------------------------------------------
%-------------------------------------------------------------------

\section{Algebraic set of atoms}

\vspace{7mm}

In this section we will concentrate on the structure of algebraic
subsets of the set $\mathcal A$ of atoms.  We will introduce the
rank of such an algebraic set and also get some connections
between \find and some usual dimensions of vector spaces over
division rings.  This will shed some new lights on these notions.

We start this section with some easy facts on algebraic sets of
atoms.  First let us recall that, in general, even for a finite
algebraic set $\De$ with basis $B$ we might have $RB_\ell \ne
R\De_\ell$ as we have seen in example \ref{algebraic example} (c).
In case of algebraic sets of atoms we have:
\begin{prop}
\label{B_l et D_l} Let $\De \subseteq \mathcal A$ be an
$F$-algebraic set with basis $B$.
\begin{enumerate}
\item[a)] $\intdel = \bigcap_{b\in B}Rb$.
\item[b)] If $\vert B \vert < \infty$ then there exist $\De_\ell$
and $B_\ell \in \mathcal F$ such that $\intdel = R\De_\ell =
\bigcap_{b\in B}Rb = RB_\ell$ and $\ell (B_\ell) = \ell (\De_\ell)
= \vert B \vert$.
\end{enumerate}
\end{prop}
\begin{proof}
a) The inclusion $\intdel \subseteq \bigcap_{b \in B}Rb$ is clear.
Now if $x\in \bigcap_{b\in B}Rb$ and $\delta \in \De$ then, thanks
to Proposition \ref{dependence of an atom} there exists a finite
subset $\Gamma$ of $B$ such that $\intgam \subseteq R\delta$ ,
hence $x\in \bigcap_{b \in B}Rb \subseteq \intgam \subseteq
R\delta$ for any $\delta \in \De$.  This shows that $\bigcap_{b
\in B}Rb \subseteq \intdel$ as desired.

b) This is clear in view of a) above and corollary \ref{finite
algebraic}.
\end{proof}

In view of the above proposition it is natural to introduce the
following notions:
\begin{defns}
\label{defn of rank}
\begin{enumerate}
\item[a)]Let $\De$ be an $F$-algebraic set of atoms and $B$ be an $F$-basis
for $\De$.  We define the {\bf rank} of $\De$, denoted $rk(\De)$,
by $rk(\De)= \vert B \vert$.
\item[b)] For $a\in R\setminus\{0\}$, let
$$V(a):=\{p\in {\mathcal {A}} | a\in Rp\}$$
\item[c)] For an algebraic subset $\De$ of $\mathcal A$ we call
the {\bf closure} of $\De$ the set of atoms which are
$F$-dependent on $\De$ and we denote this set by $\overline{\De}$.
\end{enumerate}
\end{defns}

\vspace{5mm}
\begin{lem}
\label{closure,ect} With the above notations and definitions we
have :
\begin{enumerate}
\item[a)] $V(a)$ is an $F$-algebraic set and $rk(V(a)) \leq
\ell(a)$.
\item[b)]If $\De \subseteq \mathcal A$ is an $F$-algebraic set
with an $F$-basis $B$, then $\overline \De$ is $F$-algebraic,
 $\overline B = \overline \De$
and $rk(\De) = rk(\overline \De)$. If $\De$ is of finite rank then
$\overline \De = V(\De_\ell)= V(B_\ell) = \overline B$.
\item[c)]Let $a,b \in \mathcal F \setminus \{0\}$ be such that $Ra \cap Rb
\neq 0$.  Then $V(a)\cap V(b)=\emptyset$ if and only if $Ra+Rb=R$
if and only if $\ell ([a,b]_\ell)=\ell(a)+\ell(b)$.

\item[d)] If $C$ is a finite algebraic subset of $\mathcal A$ then
$rk(C) \le \vert C \vert$ and the equality occurs if and only if
$C$ is an $F$-independent subset of $\mathcal A$.
\end{enumerate}
\end{lem}

\begin{proof}
\begin{enumerate}
\item[a)] If $a$ is a unit in $R$ then $V(a)=\emptyset$ and
so $V(a)$ is algebraic.  If $0 \ne a \in R \setminus U(R)$ we have
$0 \ne Ra\subset \bigcap_{p\in V(a)}Rp $.  This shows that $V(a)$
is an algebraic set.  Proposition \ref{B_l et D_l} b) implies that
for a finite $F$-independent set $B \subseteq V(a)$ we have $
Ra\subseteq\cap_{b\in B}Rb = RB_\ell $ and so $\vert B \vert =
\ell (B_\ell) \le \ell (a))$
\item[b)] Using Proposition \ref{dependence of an atom} it is easy to remark that
$0\neq (\cap_{\delta \in \Delta}R\delta) = (\cap_{\gamma \in
\overline\Delta}R\gamma)$ and so $\overline\Delta$ is an algebraic
set.  Obviously $\overline B \subseteq \overline\Delta$ and the
transitivity of $F$-dependence on sets of atoms gives the reverse
inclusion.  Hence $\overline B = \overline \De$ and $rk(\De) =
rk(\overline \De)$.  The last statement follows easily.

\item[c)] and d) are left to the reader.

\end{enumerate}
\end{proof}

\vspace{5mm}
\begin{thm}
\label{delta and gamma} Let $\Delta \cup \Gamma \subseteq \mathcal
A$ be an $F$-algebraic set of finite rank.  Then
\begin{enumerate}
\renewcommand{\labelenumi}{(\roman{enumi})}
\item $R(\Delta\cup \Gamma)_\ell =R\Delta_\ell \cap R\Gamma_\ell$ and
$rk(\Delta\cup \Gamma) \leq rk(\Delta)+rk(\Gamma)$.
\item $V(\Delta_\ell)\cap V(\Gamma_\ell) = \emptyset$ if and only
if equality holds in $(i)$.
\end{enumerate}
\end{thm}

\begin{proof} $(i)$ We have $R(\Delta \cup \Gamma)_\ell =
\cap_{\epsilon \in \Delta \cup \Gamma}R\epsilon = (\cap_{\delta
\in \Delta}R\delta)\bigcap (\cap_{\gamma \in \Gamma}R\gamma) =
R\Delta_\ell \cap R\Gamma_\ell$ .  The statement about rank
follows from Theorem \ref{length of gcd and llcm} and Proposition
\ref{B_l et D_l}(b).  To prove $(ii)$ assume $V(\Delta_\ell)\cap
V(\Gamma_\ell) = \emptyset$.  Then the previous lemma and $(i)$
above imply that $\ell ((\Delta \cup
\Gamma)_\ell)=\ell(\Delta_\ell)+ \ell(\Gamma_\ell)$.  Hence
$rk(\Delta\cup \Gamma) =rk(\Delta)+rk(\Gamma)$.

\end{proof}

\vspace{5mm}
\begin{thm}
Let $\Delta \cup \Gamma$ be an $F$-algebraic set in ${\mathcal
{A}}$.  Denote by $B,B'$ respectively $F$-bases for $\Delta$ and
$\Gamma$.  Then we have $\overline{\Delta} \cap \overline{\Gamma} =\emptyset$
if and only if $B \cup B'$ is an $F$-basis for $\Delta \cup
\Gamma$.
\end{thm}

\begin{proof}

  Using Lemma \ref{closure,ect} we get $\overline \De
\bigcap \overline \Gamma =\emptyset$ if and only if $\overline B
\bigcap \overline B' = \emptyset$ if and only if $\overline C
\bigcap \overline C' = \emptyset$ for any finite subsets $C, C'$
of $B$ and $B'$ respectively.  This is equivalent to $V(C_\ell)
\bigcap V(C'_\ell) = \emptyset$ i.e., using Theorem \ref{delta and
gamma}, $rk(C \cup C') = rk(C) + rk(C') = \vert C \vert + \vert C'
\vert $ for any finite subsets $C,C'$ of $B$ and $B'$
respectively.  From the above we conclude that $\overline \De
\bigcap \overline \Gamma = \emptyset $
 if and only if for any
finite subsets $C,C'$ of $B$ and $B'$ respectively we have that
$rk (C \cup C') = \vert C \cup C' \vert$ i.e. if and only if $C
\cup C'$ is $F$-independent.  Hence $\overline \De \bigcap
\overline \Gamma = \emptyset$ if and only if $B \cup B'$ is
$F$-independent.  Now, since $B$ and $B'$ are $F$-bases for $\De$
and $\Gamma$ respectively it is easy to finish the proof.
\end{proof}

  Let us recall, from section 2, that for $a,b \in R$ such
that $Ra \cap Rb \ne 0$ we wrote $Ra \cap Rb = Ra^bb = Rb^{a}a$.

\vspace{5mm}
\begin{prop}
\label{factorzation of a llcm} Let $\Gamma \subseteq \mathcal F$ be an
$F$-algebraic set of atoms in a \df\ $R$ such that $\intgam=Rh$
for some element $h \in \mathcal F$.  Then
\begin{enumerate}
\item[(i)] $h$ is a product of atoms similar to atoms in $\Gamma$.
\item[(ii)] any right atomic factor of $h$ is similar to some
atom in $\Gamma$.
\end{enumerate}
In particular, this applies to any finite subset of an $F$-
algebraic set $\De \subseteq \mathcal A$.
\end{prop}

\begin{proof} Let $B$ be a basis for $\Gamma$.  From
Corollary \ref{finite algebraic} and Proposition \ref{B_l et D_l}
we have that $\ell (h) = \vert B \vert < \infty$.  Let us put $B
=\{b_1, \dots ,b_s \}$, we will show by induction on $s$ that $h$
is a product of atoms similar to the $b_i$'s.  From \ref{B_l et
D_l} we know that $\bigcap_{i=1}^sRb_i = \intgam = Rh$.

  If $s=1$ we have $Rb_1 = Rh$ and $h$ must be an atom
associated to $b_1$.

  If $s > 1$ we have $Rh \subseteq Rb_1$ and we can write
$h = h_1b_1$.  We then have $Rh_1b_1 = Rh = Rb_1 \bigcap (\cap
_{i=2}^sRb_i)= \bigcap_{i=2}^s(Rb_1 \cap Rb_i) =
\bigcap_{i=2}^sRb_i^{b_1}b_1$.  This gives that $Rh_1 =
\bigcap_{i=2}^s Rb_i^{b_1}$.  Now $\{b_2^{b_1}, \dots ,b_s^{b_1}
\}$ is an algebraic set and the induction hypothesis implies that
$h_1$ is a product of atoms which are similar to the $b_i^{b_1}$'s
and hence similar to the $b_i$'s for $i \in \{2, \dots ,s\}$.
Since $h = h_1b_1$ we can conclude.

(ii) Let us use the same notations as in (i) above and assume that
$h=ga$ where $g \in \mathcal F$ and $a \in \mathcal A$.  We want
to show that $a$ is similar to one of the $b_i$'s.  We proceed by
induction on $s$. We have $h=h_1b_1=ga \in Ra$ with
$\bigcap_{i=2}^s Rb_i^{b_1}=Rh_1$. Hence by \ref{product formula}
either $b_1 \in Ra$ or $h_1 \in Ra^{b_1}$.  In the first case we
conclude that $a$ is associated to $b_1$ and hence $a$ and $b$ are
similar.  In the second case the induction hypothesis shows that
$a^{b_1}$ is similar to one of the $b_i^{b_1}$'s.  The
transitivity of similarity yields the conclusion.
\end{proof}

The following definition will be useful for us:

\begin{defn}
\label{full} An $F$-algebraic subset $\Gamma$ of a set $\De$ is
{\bf full} in $\De$ if any element of $\De$ which is $F$-dependent
on $\Gamma$ is already in $\Gamma$.
\end{defn}

\vspace{5mm}
\begin{lem}
\label{decomposition following a full subset}   Let $\Delta$ be an
$F$-algebraic set of atoms and $\Gamma$ be a full subset of
$\Delta$. If $Rh=\cap_{\gamma \in \Gamma}R \gamma$ and
$Rf=\cap_{\delta \in \Delta}R\delta$ then $Rf=Rgh$ where
$Rg=\cap_{d\in \Delta \setminus \Gamma}Rd^{h}$.

\end{lem}
\begin{proof}

  Since $\Gamma \subseteq \Delta$ we know that there
exists $g \in R$ such that $f=gh$, and we must show that
$Rg=\cap_{d\in \Delta \setminus \Gamma}Rd^{h}$. Now, for any $d
\in \Delta \setminus \Gamma$ we know that $f=gh \in Rd$, but since
$\Gamma$ is full in $\Delta$ we have that $h \notin Rd$ hence by
Lemma \ref{product formula} $g \in Rd^h$. This shows that $Rg
\subseteq \cap_{d\in \Delta \setminus \Gamma}Rd^{h}$.  On the
other hand, if $p \in \cap_{d\in \Delta \setminus \Gamma}Rd^{h}$
 then $ph \in \cap_{d\in \Delta \setminus \Gamma}Rd^{h}h=
 \cap_{d\in \Delta \setminus \Gamma}(Rd \cap Rh)= \cap_{d\in \Delta}Rd$
 and hence, $ph \in Rf=Rgh$. This implies that $p \in
Rg$, as required.

\end{proof}

We will study the influence of the decomposition into similarity
classes on the notions of $F$-independence and rank.   Let us
first start with the promised expression of $F$-independence of an
element in terms of the $F$-independence of the atoms appearing in
its factorization.

Let us first introduce the following notation :  for $\Delta
\subseteq R$ and $u \in R \setminus \{0\}$ we denote $\Delta^{u} =
\{g\in R \vert \exists \delta \in \Delta : Rgu = R\delta \cap Ru
\ne 0 \}$ (to justify this notation let us notice that in
\ref{a^b} we wrote $R\delta \cap Ru = R\delta^{u}u$).

\vspace{5mm}
\begin{prop}
\label{general dependence} Let $R$ be \df\ and $a=p_1p_2\cdots
p_n$ be a factorization of an element $a\in \mathcal F$ into
atoms. If $\Delta\subseteq \mathcal F$ is $F$-algebraic then $a$
is $F$-dependent on $\Delta$ if and only if either $p_n$ is
$F$-dependent on $\Delta$ or there exists $s\in\{1,2,\dots ,
n-1\}$ such that $p_s$ is $F$- dependent on $\Delta^{p_{s+1} \dots
p_n}$.
\end{prop}
\begin{proof}
Assume $a$ is $F$-dependent on $\De$.   We have $0 \ne (\intdel)
\bigcap Ra \subseteq (\intdel) \bigcap Rp_n$, so that $\De \cup
\{p_n \}$ is algebraic.  If $n=1$ the result is clear.  So let us
assume that
 $ n>1 $ and that $p_n$ is not $F$-dependent on
$\Delta$.  We leave it to the reader to check that $\Delta^{p_n}$
is algebraic.  Now there exists a finite subset $\Gamma $ of
$\Delta$ and a non unit $d \in R$ such that $R \Gamma_\ell + Ra =
Rd$.   We claim that $p_1p_2...p_{n-1}$ is $F$-dependent on
$\De^{p_n}$.  First let us remark that $\De^{p_n} \cup
\{p_1p_2...p_{n-1} \}$ is algebraic since $((\bigcap_{\delta \in
\De} R\delta^{p_n}) \bigcap Rp_1p_2...p_{n-1})p_n $
$=\intdel^{p_n}p_n \bigcap Ra = (\intdel \bigcap Rp_n) \bigcap Ra
= \intdel \bigcap Ra \ne 0$.  Now assume that $R\Gamma^{p_n}_\ell
+ Rp_1p_2...p_{n-1} = R$, then $Rp_n = (R\Gamma^{p_n}_\ell +
Rp_1p_2...p_{n-1})p_n $ $= R\Gamma^{p_n}_\ell p_n + Ra =
(R\Gamma_\ell \bigcap Rp_n) + Ra \subseteq Rd$.  Since $p_n$ is an
atom and $d$ is not a unit this leads to $Rd = Rp_n$, but then $Ra
+R\Gamma_\ell = Rp_n$ and hence $Rp_n + R\Gamma_{\ell} = Rp_n$.
This contradicts the fact that $p_n$ is not $F$-dependent over
$\Delta$ and proves the claim. Now the induction hypothesis and
the formula $(\Delta^{p_n})^q = \Delta^{qp_n}$ for any $q$ such
that $Rq \cap \intdel \ne 0 $ allow us to conclude easily.

Conversely, assume first that $p_n$ is $F$-dependent on $\De$ and
consider $\Gamma$ a finite subset of $\De$ such that $\cap_{\gamma
\in \Gamma}R\gamma + Rp_n \not = R$.  Then, $\cap_{\gamma \in
\Gamma} R\gamma + Ra \subset \cap_{\gamma \in \Gamma}R\gamma +
Rp_n \not = R$.

\noindent Now, assume that $p_n$ is $F$-dependent on $\De$, but
there exists $s \in{1,2,\dots,n-1}$ such that $p_s$ is
$F$-dependent on $\De^{p_{s+1}\cdots p_n}$.  This means that there
exists a finite subset $\Gamma_0 \subseteq \De$ such that
$\cap_{\gamma \in \Gamma_0}R\gamma^{p_{s+1}\cdots p_n}\subseteq
Rp_s$.  We want to show that $a$ is $F$-dependent on $\De$. Assume
that this is not the case.  Then, for all finite subset $\Gamma
\subseteq \De \, , \cap_{\gamma \in \Gamma} R\gamma + Ra = R$.  In
particular, $\cap_{\gamma \in \Gamma_0} R\gamma + Ra = R $.  Hence
we have $(\cap_{\gamma \in \Gamma_0} R\gamma + Ra) \cap
Rp_{s+1}\cdots p_n = Rp_{s+1}\cdots p_n$.  Since $Ra \subseteq
Rp_{s+1}\cdots p_n$, this gives $Rp_{s+1}\cdots p_n=((\cap_{\gamma
\in \Gamma_0} R\gamma) \cap Rp_{s+1}\cdots p_n) + Ra =
(\cap_{\gamma \in \Gamma_0}(R\gamma\cap Rp_{s+1}\cdots p_n)) + Ra
= \cap_{\gamma \in \Gamma_0}R\gamma^{p_{s+1}\cdots
p_n}p_{s+1}\cdots p_n + Ra$. Since $\cap_{\gamma \in
\Gamma_0}R\gamma^{p_{s+1}\cdots p_n}\subseteq Rp_s$, we finally
get $Rp_{s+1}\cdots p_n \subseteq Rp_sp_{s+1}\cdots p_n$.  This
contradiction yields the result.
\end{proof}

For an element $a$ in $R$ we denote $\Delta(a)$ the set of
elements which are similar to $a$.

\vspace{5mm}
\begin{thm}
\label{excision} Let $\Delta $ be an algebraic set of atoms in a
\df\ $R$.  If an atom $a$ is $F$-dependent on $\Delta$ then $a$ is
$F$-dependent on $\Delta \cap \Delta(a)$.
\end{thm}
\begin{proof} Since $a$ is $F$-dependent on a finite subset of $\De$
we may assume that $\De$ is finite. Put $\Gamma:=\Delta \cap
\Delta(a)$, $h:=\Gamma_\ell$ and denote by $f:=\Delta_\ell$.  We
must show that $h \in Ra$. Let us notice that for any element $d
\in \Delta \setminus \Delta(a)=\Delta \setminus \Gamma$ we have $h
\notin Rd$ (since by Proposition \ref{factorzation of a llcm} the
factors of $h$ are similar to $a$) hence $\Gamma$ is full in
$\Delta$.  Now we can write, as in the lemma \ref{decomposition
following a full subset}, $f=gh$ where $g$ is such that $Rg =
\cap_{d\in \Delta \setminus \Gamma}Rd^{h}$.  For any $d \in \Delta
\setminus \Gamma$, we must have $g \in Rd^h$.  On the other hand,
since $a$ is $F$-dependent on $\Delta$, we have that $gh\in Ra$.
Assume now that $h\notin Ra$.  Then, thanks to Lemma \ref{product
formula}, $g \in Ra^h$ . But this would mean that an element of
$\Delta(a)$ is a factor of $g$. This contradicts the definition of
$g$ and shows that $h$ must be in $Ra$, as desired.

\end{proof}
This theorem has an immediate useful corollary which will
essentially reduce the study of an algebraic set of atoms to the
case of an algebraic set contained in a similarity class.

\vspace{5mm}
\begin{cor}
\label{decomposition into similarity classes} Let $\Delta \subset
\mathcal A$ be an algebraic set of finite rank.  Then $\De$
 intersects a finite number of similarity classes.  More precisely
 : if $r=rk(\De)$, there exist $n\le r$ non similar atoms $
 p_1,\dots ,p_n \in \mathcal A$ such that $\De =
 \bigcup_{i=1}^{n} \De_i $ where $\De_i= \De \cap \De(p_i)$ for
 $i \in \{1, \dots ,n \}$.  Moreover if $B_i$ is an $F$-basis for $\De_i$
 then $B:=\bigcup B_i$ is an $F$-basis for $\De$ and
 $$
rk(\De) = \sum_{i=1}^{n} rk(\De_i)  \quad \quad ; \quad \quad
\overline{\De} = \bigcup_{i=1}^n \overline{\De_i}
 $$
 In particular, if $f \in \mathcal F$
 then $V(f)$ intersects at most $l(f)$ similarity classes.
 \end{cor}
\begin{proof}
Assume at the contrary that $\De$ intersects more than $r=rk(\De)$
similarity classes and let $a_1,\dots,a_{r+1}$ be elements of
$\De$ belonging to distinct similarity classes.  Then the above
theorem shows that $\{a_1,\dots,a_{r+1}\}$ are $F$-independent, hence
$rk(\De) \ge r+1$, a contradiction.  Now if $x \in \De$, then the
above theorem shows that $x$ is $F$-dependent on $\De_i$ for some
$i\in \{1,\dots,n\}$; i.e. $x$ is $F$-dependent on some $B_i$. On
the other hand if $y\in \bigcup_{i=1}^{n} B_i$ is $F$-dependent on
$\bigcup_{i=1}^{n} B_i \setminus \{y\}$ then $y \in B_i$ for some
$i$ and is $F$-dependent on $B_i \setminus \{y\}$.  This
contradiction allows us to conclude that $B$ is an $F$-basis for
$\De$.

  It remains to prove that $\overline{\De} \subseteq
\bigcup_{i=1}^n \overline{\De_i}$ (the other inclusion being
obvious).  Let $p \in \mathcal A$ be an element which is
$F$-dependent on $\De = \bigcup_{i=1}^n \De_i$.  By the above
theorem \ref{excision} we know that $p$ is $F$-dependent on
$\De(p) \cap (\bigcup_{i=1}^n \De_i) = \bigcup_{i=1}^n (\De_i \cap
\De(p)) $.  Since the $p_i$'s are non similar all but one of these
intersections are empty and so there exists $j\in\{1,\dots,n\}$
such that $p$ is $F$-dependent on $\De_j$.

\end{proof}

  The notion of $F$-independence
will be particularly explicit inside the similarity classes
$\Delta(p_i)$.  Let us recall that for an atom $p \in \mathcal A$,
the ring $End_R(R/Rp)$, denoted $C(p)$, is in fact a division ring
(Cf. Corollary \ref{Schur}).  It turns out that in the similarity
class $\Delta(p)$ of an atom $p$ the notion of $F$-independence
can be translated in terms of usual linear dependence over this
division ring $C(p)$.  Let us also recall that $R/Rp$ has a
natural structure of right $C(p)$-vector space.  In the following
definition we introduce a very useful map.
\begin{defn}
\label{lambda map} Let $p$ be an atom and $f \in R$.  We define
$$
\lambda_{f,p}:R/Rp \longrightarrow R/Rp : x + Rp \mapsto fx + Rp
\, .
$$
\end{defn}

\vspace{5mm}
\begin{thm}
\label{F-independence and independence over End(R/Rp)}
\begin{enumerate}
\item [(a)] Let $f$ be an element in $\mathcal F$ the map
$\lambda_{f,p}$ is a right $C(p)$-linear map with
$Ker(\lambda_{f,p}) = \{x + Rp \;\vert \; fx \in Rp\}$ and we have
$$dim_{C(p)}Ker(\lambda_{f,p}) \le \ell(f).$$
\item [(b)] Let $\Delta$ be an algebraic subset contained in the similarity
class $\Delta(p)$ of an atom $p$.  Let $\{p_1,p_2,...,p_n \}$ be a
subset of $\Delta$ and for $i\in \{1,2,...n\}$ let $\phi_i :
R/Rp_i \longrightarrow R/Rp$ be isomorphisms of left $R$-modules.
Then the set $\{p_1,p_2,...,p_n\}$ is $F$-dependent if and only if
the set $\{\phi_1(1+Rp_1),\phi_2(1+Rp_2),...,\phi_n(1+Rp_n)\}$ is
right $C(p)$-dependent.
\item[(c)] For $f \in \mathcal F$ and $p \in \mathcal A$, we have
$$dim_{C(p)}Ker(\lambda_{f,p}) = rk(V(f) \cap \Delta (p)). $$
\end{enumerate}
\end{thm}
\begin{proof}
 (a) We leave it to the reader to check that for the
natural structure of right $C(p)$-vector space on $R/Rp$, the map
$\lambda_{f,p}$ is a right homomorphism.  The given description of
$ker \lambda_{f,p}$ is straightforward and we only need prove that
$dim_{C(p)}ker\lambda_{f,p} \le \ell(f)$.  We proceed by induction
on $\ell(f)$.  The claim is obvious if $\ell(f)=0$.  If $f$ is an
atom and $x+Rp,\; y+Rp$ are nonzero elements in $ker
\lambda_{f,p}$ then $fx \in Rp$ and $fy \in Rp$.  Using the
notations of Lemma \ref{product formula} we have $f \in Rp^x \cap
Rp^y$.  Since $x \notin Rp$ and $y \notin Rp$, $p^x$ and $p^y$ are
not units in $R$ and, $f$ being an atom we conclude that $f=\alpha
p^x =\beta p^y$ for units $\alpha$ and $\beta$ in $R$.  Define the
isomorphisms $\phi_x : \frac{R}{Rp^x} \longrightarrow
\frac{R}{Rp}: 1+Rp^x \mapsto x + Rp$ and similarly for $\phi_y$.
Now, since $Rp^x=Rf=Rp^y$ we have that
$\frac{R}{Rp^x}=\frac{R}{Rp^y}$ and the map $\gamma
:=(\phi_x)^{-1} \circ \phi_y \in End_R(R/Rp)$ is such that $\gamma
(x+Rp) = y+Rp$. Hence $x +Rp$ and $y + Rp$ are right
$C(p)$-dependent. This shows that $dim_{C(p)}ker \lambda_{f,p} \le
1$ as desired. For the general case we remark that if $f=f_1f_2...
f_r$ is an atomic decomposition of $f$ then we have
$\lambda_{f,p}=\lambda_{f_1,p} \circ \lambda_{f_2,p}\dots \circ
\lambda_{f_r,p}$.  Hence
 $dim_{C(p)} ker \lambda_{f,p} \le \sum_{i=1}^{r}
 dim_{C(p)}ker\lambda_{f_i,p} \le r = \ell(f)$, as desired.

(b) Let us put $x_i:=\phi_i(1+Rp_i)$.  We then have $p_ix_i = 0
\in R/Rp$.  First let us assume that the $x_i$'s are right
$C(p)$-dependent and let $\sum_{i=1}^{n} x_i\gamma_i =0 $ be a
dependence relation.  Without loss of generality we may assume
that $\gamma_n \ne 0$ and thus write $x_n =
\sum_{i=1}^{n-1}x_i\psi_i$ for some $\psi_i\in C(p)$.  Since
$\Delta$ is algebraic there exists $f$ in $R$  such that
$\bigcap_{i=1}^{i=n-1}Rp_i = Rf$ and we will show that $Rf
\subseteq Rp_n$.  We know there exist $f_1,f_2, \dots,f_{n-1}$
such that $f=f_ip_i$ for $i=1,2, \dots,n-1$ and since $p_ix_i = 0$
in $R/Rp$, we get $fx_i=f_ip_ix_i =0$ for $i=1,2, \dots,n-1$.
This leads to $\phi_n(f + Rp_n) = fx_n =
\sum_{i=1}^{n-1}fx_i\psi_i = 0$ and so $f \in Rp_n$, as desired.

  Conversely let us suppose that $p_1,p_2, \dots,p_n$ are
$F$-dependent.  Since these elements are contained in an algebraic
set, we have $\cap_{i=1}^{n}Rp_i = Rf $ for some $f\in \mathcal F$
and since they are $F$-dependent we know by proposition \ref{B_l
et D_l}(b) that $\ell(f) \le n-1$.  For $i=1,2, \dots ,n$ let us
write $f = f_ip_i$.  Now, since for $i=1,2, \dots ,n$ we have
$p_ix_i = 0 \in R/Rp$, we have $\lambda_{f,p}(x_i) = fx_i =
f_ip_ix_i =0 $ i.e. $\{x_1, \dots,x_n \} \subseteq
ker\lambda_{f,p}$.  By part (a) above we have that
$dim_{C(p)}ker\lambda_{f,p} \leq \ell(f) \leq n-1$ and we conclude
that $\{x_1,x_2,\dots ,x_n \}$ are right $C(p)$-dependent.

(c) Let $\{p_1,\dots,p_n\}$ be an $F$-basis for $V(f) \cap
\Delta(p)$ and put $R/Rp_i \stackrel{\phi_i}{\cong} R/Rp : y+
Rp_i \mapsto yx_i + Rp $ for some $x_i\in R$.  Let $y_i \in R$
be such that $p_ix_i = y_ip$; since $f \in Rp_i$, we have $fx_i
\in Ry_ip$ and so $\phi_i(1+Rp_i) = x_i +Rp\in ker(\lambda_{f,p})$
and part (b) above shows that these elements are
$C(p)$-independent.  We thus conclude that $rk(V(f)\cap\Delta(p))
\le dim_{C(p)}ker(\lambda_{f,p})$.

  Conversely if $x_1+Rp,\dots,x_n+Rp$ are
$C(p)$-independent in $ker(\lambda_{f,p})$ then $fx_i \in Rp$ and
since $x_i \notin Rp$ we get $f \in Rp^{x_i}$ and from part b)
again we easily conclude that $p^{x_1},\dots,p^{x_n}$ are
$F$-independent elements in $V(f) \cap \Delta(p)$.
\end{proof}

Part a) in the above theorem was obtained by P.M.Cohn
\cite[Theorem 5.8, P.233]{C1} and part b) was inspired by similar
results obtained for Ore extensions \cite{LL1}.

  With the help of the previous theorem we are ready to
present, as a corollary, the full computation of the rank of an
algebraic subset $\De \subseteq \mathcal A $ as well as the
description of the closure $\overline{\De}$. Recall that for an
algebraic set of finite rank corollary \ref{decomposition into
similarity classes} shows that $\De$ intersects a finite number of
similarity classes $\De(p_1), \dots,\De(p_n)$ and we can write
$\De =
 \bigcup_{i=1}^{n} \De_i $ where $\De_i= \De \cap \De(p_i)$ for
 $i \in \{1, \dots ,n \}$.  Now, for any $\gamma \in \De_i$ let
 $\phi_\gamma$ be an isomorphism $R/R\gamma \cong R/Rp_i$ and
 denote by $Y_i$ the right $C(p_i)-$subspace of $R/Rp_i$ defined by
 $Y_i := \sum_{\gamma \in \De_i} \phi_\gamma(1+R\gamma)C(p_i)$.
 With these notations we can state:

\vspace{5mm}
\begin{cor}
\label{rank decomposition} Let $\De$ be an algebraic subset of
$\mathcal A$.   Then $rk(\De) = \infty$ if and only if one of the
following holds :

a)$\De$ contains infinitely many non similar atoms.

b)There exist $p\in \mathcal A$ and infinitely many atoms $p_i$ in
$\De \cap \De(p)$ such that their images into $R/Rp$ by the
isomorphisms $\phi_i :R/Rp_i \cong R/Rp$ generate an infinite
dimensional vector space over $C(p)$.

If none of these conditions is satisfied then $\De$ is of finite
rank and, using the above notations,  we have :

 $$ \De = \bigcup_{i=1}^{n} \De_i \quad  \quad rk(\De)=
\sum_{i=1}^{n}dim_{C(p_i)}Y_i  \quad \quad \overline{\De} =
\bigcup_{i=1}^{n}\overline{\De_i} .$$

In particular if $f \in \mathcal F$ and $V(f)= \cup_{i=1}^r (V(f)
\cap \Delta(p_i))$ is the decomposition of $V(f)$ into similarity
classes one has

$$ rkV(f)=\sum_{i=1}^rdim_{C(p_i)}ker(\lambda_{f,p_i}).$$
\end{cor}
\begin{proof}
The proof uses \ref{decomposition into similarity classes} and
theorem \ref{F-independence and independence over
End(R/Rp)}(b),(c).
\end{proof}

\vspace{5mm}
The next result, although a bit technical, will be
helpful.

\begin{prop}
\label{bases and independence} Let $h$ be a nonzero element in $R$
and $\{a_1,\ldots ,a_n \}$ be an $F$-basis for $V(h)$.  An
algebraic set $\{a_1,\dots ,a_n,b_1, \dots ,b_m \} \subseteq
\mathcal A$ such that $\{b_1, \dots ,b_m\} \subset R \setminus
V(h) $ is $F$-independent if and only if $\{b^{h}_1,\dots,b^{h}_m
\}$ is $F$-independent.
\end{prop}
\begin{proof}
If $h$ is a unit $ V(h) $ is empty and $\{b_1,b_2,\dots,b_n \}$ is
$F$-independent if and only if $ \{b_1^h,b_2^h, \dots,b_n^h \} $
is $F$-independent.  We may thus assume that $h$ is not a unit and
we begin with the `` only if'' part.  By Theorem \ref{excision},
we know that elements in different similarity classes are
$F$-independent.  Lemma \ref{a^b} shows that $b^{h}_j \sim b_j$,
hence we may assume that $\{a_1,\dots,a_n,b_1,\dots,b_m \}$ is
contained in a single conjugacy class, say $\Delta(p)$.  Let, for
$i\in \{ 1,\dots,n \}$ and $j\in \{1,\dots,m \}$,
\begin{itemize}
\item $\phi_i : R/Ra_i \longrightarrow R/Rp :1+Ra_i\mapsto \alpha_i +Rp,$
\item $\psi_j : R/Rb_j \longrightarrow R/Rp :1+Rb_j\mapsto \beta_j +Rp$
\item $\sigma_j : R/Rb^{h}_j \longrightarrow R/Rb_j :1+Rb^{h}_j \mapsto h+Rb_j$
\end{itemize}
be isomorphisms of left $R$-modules.
%By Theorem
%\ref{F-independence and independence over End(R/Rp)}
%(b), we know that $\{\alpha_1 +Rp,\dots,\alpha_n +Rp,\beta_1
%+Rp,\dots,\beta_m +Rp\}$ is right $C(p)$-independent.  For $j\in
%\{1,\ldots,m \}$ define the left $R$-isomorphisms $\sigma_j :
%R/Rb^{h}_j \longrightarrow R/Rb_j \; : 1+Rb^{h}_j \mapsto h+Rb_j$.
Then $\psi_j \circ \sigma_j$ is an $R$-isomorphism of left modules
between $R/Rb^{h}_j$ and $R/Rp$ such that $\psi_j \circ
\sigma_j(1+Rb^{h}_j)=h\beta_j + Rp$.  Now, assume that
$\{b^{h}_1,\ldots,b^{h}_m \}$ is $F$-dependent.  Then Theorem
\ref{F-independence and independence over End(R/Rp)} (b) shows
that $\{h\beta_1 + Rp,\ldots, h\beta_m + Rp \}$ is right
$C(p)$-dependent.  So, there exist $\eta_1 ,\ldots,\eta_m \in
C(p)$ not all zero such  that $$Rp=\sum_{j=1}^m ((h\beta_j
+Rp)\eta_j)=h(\sum_{j=1}^m(\beta_j +Rp)\eta_j).\label{*}$$  Let us write
$(\beta_j +Rp)\eta_j = \beta^{'}_j +Rp$.  So we get
$$h(\sum_{j=1}^m
\beta^{'}_j)\in Rp.$$  Let us remark that by Theorem
\ref{F-independence and independence over End(R/Rp)} (b), we know
that $\{\beta_1 +Rp,\dots,\beta_m +Rp\}$ is right
$C(p)$-independent. This implies that $\sum_{j=1}^m \beta^{'}_j
+Rp \neq Rp$, and so, $\sum_{j=1}^m \beta^{'}_j \notin Rp$.  Using
Lemma \ref{product formula}, we get $h\in Rp^{\sum \beta^{'}_j}$.
This shows that $p^{\sum \beta^{'}_j} \in V(h)$ and so,
$\{p^{\sum\beta^{'}_j},a_1,\ldots,a_n \}$ is $F$-dependent.
Considering the $\phi_i$'s and the isomorphism of left $R$-modules
$R/Rp^{\sum \beta^{'}_j} \longrightarrow R/Rp : 1+Rp^{\sum
\beta^{'}_j} \mapsto \sum_{j=1}^m \beta^{'}_j +Rp$ , Theorem \ref
{F-independence and independence over End(R/Rp)} (b) again shows
that $\{\sum_{j=1}^m \beta^{'}_j +Rp,\alpha_1 +Rp,\ldots,\alpha_n
+Rp \}$ is right $C(p)$-dependent.  In other words $\sum_{j=1}^m
\beta^{'}_j +Rp=\sum_{j=1}^m (\beta_j +Rp)\eta_j$ is right
$C(p)-$dependent on $\{\alpha_1 +Rp,\ldots,\alpha_n +Rp \}$ and so
$\{\alpha_1 +Rp,\ldots,\alpha_n +Rp,\beta_1 +Rp,\ldots,\beta_m +Rp
\}$ is also right $C(p)$-dependent.  This gives a contradiction,
by Theorem \ref{F-independence and independence over End(R/Rp)}
(b), since $\{a_1,\dots,a_n,b_1,\dots,b_m \}$ is $F$-independent.

For the ``if'' part assume $\{b^{h}_1,\ldots,b^{h}_m \}$ is
$F$-independent but $\Delta=\{a_1,\dots,a_n,$ $b_1,\ldots,b_m \}$
is $F$-dependent.  Let us suppose that $a_i$ is $F$-dependent on
$\Delta\setminus \{a_i \}$.  Let $\Delta_i$ be a minimal subset of
$\Delta\setminus \{a_i \}$ such that $a_i$ is $F$-dependent on
$\Delta_i$.  As $\{a_1,\ldots,a_n \}$ is $F$-independent, some
$b_j$ belongs to $\Delta_i$.  Now Proposition \ref{exchange} shows
that $b_j$ is $F$-dependent on
$(\Delta_i\cup\{a_i\})\setminus\{b_j\}$.  So we may assume that
some $b_j$ is $F$-dependent on $\{a_1,\dots,a_n,b_1,\ldots,b_m
\}\setminus \{b_j \}$, say $b_m$.  Let us define
$Rf:=R[b^{h}_1,\ldots,b^{h}_{m-1}]_\ell h$.  Thanks to Lemma
\ref{product formula}, we know that $f$ is a least left common
multiple of $\{a_1,\dots,a_n,b_1,\ldots,b_{m-1} \}$.  As $b_m$ is
$F$-dependent on $\{a_1,\dots,a_n,b_1, \ldots,b_{m-1} \}$, $f$ is
also a left multiple of $b_m$.  But, since $b_m \notin V(h)$,
Lemma \ref{product formula} shows that
$R[b^{h}_1,\ldots,b^{h}_{m-1}]_\ell \subset Rb^{h}_m$ i.e.
$b^{h}_m$ is $F$-dependent on $\{b^{h}_1,\ldots,b^{h}_{m-1} \}$.
This gives a contradiction.
\end{proof}

\vspace{10mm}

\section{fully reducible elements}
%----------------------------------------------------------------------
%----------------------------------------------------------------------
\vspace{7mm}

\begin{defn}
\label{fully reducible} An element $f\in \mathcal F$ is fully
reducible if there exist atoms $p_1,\dots,p_n\in R$ such that $Rf
= \bigcap_{i=1}^{n} Rp_i$
\end{defn}

This notion was introduced by Ore for skew polynomials \cite{Ore}
and for \dfs\ by P.M.Cohn \cite{C1}. It was also used for
product of linear polynomials in Ore extensions (under the name of
separate zeros) by J.Treur \cite{T} and G.Cauchon \cite{Ca} and
(under the name of Wedderburn polynomials) by T.Y.Lam and A.Leroy
\cite{LL1} and \cite{LL2}.

%***in fact in Cohn's book this notion is defined by taking any
%intersection of $Rp_i$'s *not necessarily a finite one.  In fact
%if $f$ is fully reducible then $f$ is of finite length and more
%precise results are given in the next theorem. Maybe we could just
%assume $f\in \mathcal F$.   of course for an atom $p$ , $R/Rp$
%need not be simple (give an example) but if we can characterize
%when this occurs we could relate some fully reducible elements to
%semisimplicity and the density theorem would give some interesting
%conclusions.  On the other hand maybe there is a kind of more
%"general" density theorem that would work here...If $f$ is fully
%reducible say $Rf = \bigcap Rp_i$ what can be said about $R/Rf$ ?
%Notice that if $R/Rf$ is semisimple then the density thm. would
%imply that    what is $\De(f)$ in terms of $\De (p_i)$ when are
%the $\De(p_i)$ algebraic ?  When is $\De(f)$ algebraic ($f$ fully
%irreducible) ? ***
\vspace{5mm}

 The set of fully reducible elements will be denoted
by $\mathcal R$.

\begin{lem} \label{injection,surjection,fully reducible} Let
$f,g$ be nonzero elements of a \df\ $R$ and suppose
that $g \in \mathcal R$.  Then
\begin{enumerate}
\item[a)] If $\phi : R/Rf \longrightarrow R/Rg$ is an injective
$R$-morphism then $f \in \mathcal R$.
\item[b)] If $\psi : R/gR \longrightarrow R/fR$ is a surjective
$R$-morphism then $f \in \mathcal R$.
\end{enumerate}
  In particular, if $f \sim g$ then $f \in \mathcal R$ and
in this case if $Rg = \cap_{i=1}^{n} Rp_i$, then $Rf =
\cap_{i=1}^{n} Rp'_i$ where, for $1\le i \le n \, , \; p'_i \sim
p_i$.

\end{lem}
\begin{proof}
a) Let $x \in R$ be such that $\phi(1 + Rf)= x + Rg$ and let $y
\in R$ be such that $fx=yg$.  Lemma \ref{injectivity} shows that
$\phi$ is injective if and only if $Rx \cap Rg = Rfx$. Since, by
hypothesis, $g \in \mathcal R$ there exist atoms $p_i$'s such that
$Rg = \cap_{i=1}^nRp_i$.  We thus have $Rfx =Rx \cap (\cap_{i=1}^n
Rp_i) = \cap_i(Rx \cap Rp_i) = \cap_i Rp_i^xx$. Hence we get $Rf =
\cap_i Rp_i^x$.  This yields that $f$ is fully reducible, as
requested.

b) This follows from Lemma \ref{injectivity}.

The particular case is due to the fact that in the above proof
$p_i^x\sim p_i$.
\end{proof}

Before stating the next theorem let us mention a nice consequence
of the above lemma based on the results of section 1. \vspace{5mm}
\begin{cor}
\label{g fully implies g^f fully } Let $f,g$ be nonzero
elements of a \df\ $R$ and suppose that $g \in \mathcal
R$.  If $Rf \cap Rg = Rg'f$ then $g' \in \mathcal R$.

In particular, with our standard notation, we have $g^f \in
\mathcal R$.
\end{cor}
\begin{proof}
This is an easy consequence of Lemmas
\ref{injection,surjection,fully reducible} and \ref{injectivity}.

The particular case is merely a translation of the statement using our
previous notation.
\end{proof}

Let us now come to the promised theorem showing that the notion of
reducibility is symmetric.  A constructive proof was given in
\cite[Theorem 3.6]{LL}.  We include here a short one based on
Lemma \ref{injection,surjection,fully reducible}.

\vspace{5mm}
\begin{thm}
\label{left-right symmetry}
 Suppose $\,0\neq Rf=Rp_1\cap
\dots \cap Rp_n\,$ is an irredundant intersection, where the
$\,p_i\,$'s are atoms in $\,R$.  If we write $\,\bigcap_{j\neq
i}\,Rp_j=Rg_i\,$ and $\,f=p'_ig_i\,$ $(1\leq i\leq n)$, then
\begin{enumerate}
\item[(1)] for each $\,i$, $\,p'_i\,$ is an atom similar to
$\,p_i\,$;
\item[(2)] $\,fR=\bigcap^n_{i=1}\,p'_iR\,$;
\item[(3)] the intersection representation for $\,fR\,$ in (2)
is irredundant.
\end{enumerate}
\end{thm}
\begin{proof}
We will proceed by induction on $n$.  If $n = 1$, $f=up_1=p_1' \,
,\; u\in U(R)\, , $ and $fR = p_1'R$.  Now, if $n>1$, Lemma
\ref{similarity} shows that $Rf=Rp_1 \cap Rp_2 \cap \dots \cap
Rp_n = Rp_1 \cap Rg_1 = Rp'_1g_1 =Rg'_1p_1 \; ,\; p_1 \sim p'_1 \;
, \; g_1 \sim g'_1$ and $fR = g'_1R \cap p'_1R$. Since $Rg_1 =
\cap_{j\geq 2}Rp_j$, we know that $g_1$ is fully reducible and the
above lemma \ref{injection,surjection,fully reducible} shows that
$g'_1$ is also fully reducible i.e. $Rg'_1 = \cap_{j\geq 2}Rq_j$
where the $q_j$'s are similar to the $p_j$'s.  The induction
hypothesis then gives $g'_1R = \cap_{j\geq 2}p'_jR$ where $p'_j$
are atoms and $p'_j \sim q_j \sim p_j$.  We then get $fR = g'_1R
\cap p'_1R= \cap_{i=1}^{n}p'_iR$, with $p_i \sim p_i' $ for $1\le
i \le n$, as desired.

\end{proof}
\vspace{5mm}
\begin{cor} \label{injection,surjection;fully reducible} Let
$f,g$ be nonzero elements of a \df\ $R$ and suppose
that $g \in \mathcal R$. Then
\begin{enumerate}
\item[a)] If $\phi : R/fR \longrightarrow R/gR$ is a an injective
$R$-morphism then $f \in \mathcal R$.
\item[b)] If $\psi : R/Rg \longrightarrow R/Rf$ is a surjective
$R$-morphism then $f \in \mathcal R$.
\end{enumerate}

\end{cor}
\vspace{5mm}
\begin{cor}
\label{g fully implies ^gf fully } Let $f,g$ be nonzero elements
of a \df\ $R$ and suppose that $g \in \mathcal R$.  If $fR \cap gR
= fg'R$ then $g' \in \mathcal R$.
\end{cor}
\vspace{5mm}

The following result is easy but useful :

\begin{lem}
\label{F-basis for v(f)} Let $\{p_1,\dots,p_n\} \subseteq \mathcal
A$ be a finite set of atoms and $f$ an element of $ \mathcal F$.
The following are equivalent :
\begin{enumerate}
\item[i)] $Rf = \cap_{i=1}^nRp_i$ where the intersection is irredundant.
\item[ii)] $n=l(f)$ and $\{p_1,\dots,p_n\}$ is an $F$-basis
for $V(f)$.
\end{enumerate}
In particular, $f\in R$ is fully reducible if and only if
$rkV(f)=\ell(f)$.
\end{lem}
\begin{proof}
i) $\Rightarrow$ ii) Of course, $p_i \in V(f)$ and if $a\in V(f)$
then $\cap_{i=1}^nRp_i = Rf \subseteq Ra$ and Proposition
\ref{dependence of an atom} shows that a is $F$-dependent on
$B:=\{p_1,\dots,p_n\}$.  This means that $V(f)$ is $F$-dependent
on $B$.  The fact that the intersection is irredundant implies
that $B$ is an $F$-independent subset of $R$, and the conclusion
follows.

  ii) $\Rightarrow$ i) Obviously we have $Rf \subseteq
\cap_{i=1}^nRp_i$ and this last intersection is irrredundant since
the set $\{p_1,\dots,p_n\}$ is $F$-independent.  There exists
$g\in R$ such that $\cap_{i=1}^nRp_i = Rg$.  The implication
proved above shows that $l(g) = n = l(f)$ and we conclude that $Rg
= Rf$.

  The final statement is now obvious.
\end{proof}

 In the next theorem we will give a few more
characterizations of fully reducible elements and further analyze
the structure of the set ${\mathcal{R}}$ of these elements.  In
this theorem we will use the following notations: $\Delta (p)$
will stand for the similarity class determined by an element $p$.
For an element $f\in R$ we will write as in \ref{decomposition
into similarity classes} and \ref{rank decomposition}\ $V(f) =
\cup_{i=1}^{r} \Delta_i $ where for $i\in \{1,2,\dots ,r \}$, $
\Delta_i = V(f) \cap \Delta (q_i)$ is the intersection of $V(f)$
with the similarity class $\Delta(q_i)$ of some atoms $q_i \in
\mathcal A$.  By the term a "factor" of $f\in R$ we mean an
element $g\in R\setminus U(R)$ such that there exist $p,q \in R$
with $f=pgq$.  We say that $g$ and $h$ are neighbouring factors of
an element $f$ if there exist $p,q \in R$ such that $f= pghq$.
Let us recall from Corollary \ref{Schur} that for any $i\in
\{1,2,...,r \} $, $C(q_i) := End(R/Rq_i)$ is a division ring and
remark that $R/Rq_i$ is a right $C(q_i)$-vector space.

%\begin{lem}
%Let $a,b$ be nonzero atoms in $R$. Then
%$dim_{C(b)}Ker(\lambda_{a,b})$ is equal to $1$ if $a \sim b$ and
%equal to $0$ if not.
%\end{lem}

\vspace{5mm}
\begin{thm}
\label{characterizations of fully reducible} Let $R$ be a {\df}
and let $f\in \mathcal F$.  Then the following are equivalent:
\begin{enumerate}
\item[(i)] $f$ is fully reducible.
\item[(ii)] $rkV(f) = \ell(f)$.
\item[(iii)] Let $V(f)= \cup_{i=1}^r(V(f) \cap \Delta(q_i)))$
be the decomposition of $V(f)$ into similarity classes then $\ell
(f) = \sum_{i=1}^r dim_{C(q_i)}ker(\lambda_{f,q_i})$, \\ where
$C(q_i)=End_R(R/Rq_i)$ is a division ring.
\item[(iv)]  There exist atoms $p_1,p_2,\dots,p_n$ such that
$R/Rf \cong \bigoplus_{i=1}^nR/Rp_i$.
\item[(v)] All factors of $f$ are fully reducible.
\item[(vi)] Every product of two neighbouring factors of $f$ is
fully reducible.
\item[(vii)] Every product of two neighbouring atomic factors of $f$ is
fully reducible.
\item[(viii)]For any $g \in R$ if  $V(f) \subseteq V(g)$ then $g\in
Rf$.
\end{enumerate}
\end{thm}
\begin{proof}
$(i)\Leftrightarrow (ii)$ This comes from Lemma \ref{F-basis for
v(f)}.

  $(ii) \Leftrightarrow (iii)$ This is an immediate
consequence of Corollary \ref{rank decomposition}.

  $(i)\Longrightarrow (iv)$ Assume $Rf = \cap_{i=1}^nRp_i$
where $p_i \in \mathcal A$ and the intersection is irredundant. We
shall show, by induction on $n$, that $R/Rf \cong
\oplus_{i=1}^{n}R/Rp_i$.  If $n=1$, the result is clear.  Let us
write $Rf_n = \cap_{i=1}^{n-1}Rp_i$.  We then have $Rf = Rp_n \cap
Rf_n$ and $Rp_n + Rf_n = R$ so that $R/Rf \cong R/Rp_n \oplus
R/Rf_n$.  The induction hypothesis gives $R/Rf_n \cong
\bigoplus_{i=1}^{n-1}R/Rp_i$ and enables us to conclude.

  $(iv) \Longrightarrow (i)$ We assume $R/Rf
\stackrel{\phi}{\longrightarrow} \bigoplus_{i=1}^nR/Rp_i$ is an
isomorphism.  Let $x_1,\dots,x_n \in R$ be such that $\phi(1 + Rf)
= (x_1 + Rp_1,\dots,x_n + Rp_n)$.  Since $\phi$ is well defined
and onto we have, for all $i\in \{1,\dots,n\}, \; fx_i \in Rp_i$
and $x_i \notin Rp_i$.  This leads to the fact that $f\in
\cap_{i=1}^{n} Rp_i^{x_i}$.  Hence there exists a $g\in R$ such
that $Rf \subseteq \cap_{i=1}^{n} Rp_i^{x_i} = Rg$.  In particular
we have $gx_i \in Rp_i$ for all $i\in \{1,\dots,n\}$ and $\phi (g
+ Rf) = 0$. Since $\phi$ is injective we conclude that $g\in Rf$
and $Rg \subseteq Rf$.  This shows that $Rf = Rg= \cap_{i=1}^{n}
Rp_i^{x_i}$ and, since the $p_i^{x_i}$'s are atoms, we have that
$f \in \mathcal R$, as desired.

  $(i) \Rightarrow (v)$ Assume $f=gh$.  We then have an
injective map of left $R$-modules : $R/Rg
\stackrel{.h}{\longrightarrow} R/Rf$ and Lemma
\ref{injection,surjection,fully reducible} shows that $g \in
\mathcal R$.  Similarly the injective map of right $R$-modules
$R/hR \stackrel{g.}{\longrightarrow} R/fR$ implies that $h \in
\mathcal R$.  The case of a middle factor is then clear.

  $(v) \Rightarrow (vi)$ and $(vi) \Rightarrow (vii)$
These are clear.

  $(vii) \Rightarrow (ii)$  We proceed by induction on
$n=\ell(f)$.   If $n=1$, $f$ is an atom hence belongs to $\mathcal
R$.  If $n>1$ we can write $f=ga$ for some $a \in \mathcal A$ and
$g\in R$ such that $\ell(g)=n-1$.  Clearly $g$ also satisfies the
condition in $(vii)$ and the induction hypothesis implies that $g
\in \mathcal R$.  Let us write $Rg = \cap_{i=1}^{n-1}Rp_i$ where
the $p_i$'s are in $\mathcal A$ and form an $F$-basis for $V(g)$
(cf. lemma \ref{F-basis for v(f)}).  Then $Rga =
\cap_{i=1}^{n-1}Rp_ia$ and the hypothesis shows that $p_ia \in
\mathcal R$ so that there exist $c_1,\dots,c_{n-1} \in R$ with
$Rp_ia = Rc_i \cap Ra =Rc_i^{a}a$ and we get $Rp_i=Rc_i^{a}$ which
shows that the $c_i^{a}$'s form an $F$-basis for $V(g)$.
Proposition \ref{bases and independence} then implies that
$\{c_1,\dots,c_{n-1},a\}$ are $F$-independent.  Remarking also
that $\{c_1,\dots,c_{n-1},a\} \subseteq V(ga)=V(f)$, we thus have
$rk(V(f)) \ge n = \ell(f)$.  Since the inequality $rk(V(f)) \le
\ell(f)$ is always true we get that $rk(V(f))=\ell (f)$, as
desired.

  $(i) \Rightarrow (viii)$  Assume $f \in \mathcal R$ and
let us write $Rf = \cap_{i=1}^n Rp_i$.  Hence $p_i \in V(f)
\subseteq V(g)$ and $g \in \cap_{i=1}^nRp_i=Rf$.

  $(viii) \Rightarrow (ii)$  Let us put $\ell(f) = n$.  By
$(viii)$ any element which is a left common multiple of an
$F$-basis of $V(f)$ has length $\ge n$ thus $rkV(f)\ge n$.  Since
the converse inequality always holds we get $(ii)$.

\end{proof}

\vspace{5mm}
\begin{rem}
\label{remark on thm.} It is worth to mention the relations
between the $p_i$'s and the $q_i$'s appearing in the above
theorem.  First let us notice that it is clear from the proof
that, if $f$ is fully reducible and $Rf = \cap_{i=1}^nRp_i$ is an
irredundant representation where the $p_i$'s are atoms, then these
atoms are exactly those appearing in statement $(iv)$ of the
theorem.  Let us also recall that we know from \ref{F-basis for
v(f)} that these atoms form an $F$-basis for $V(f)$.   It is then
clear that every similarity class intersecting non trivially
$V(f)$ contains at least one of the $p_i$'s.   Since the $q_i$'s
must represent these similarity classes we can just choose the
$q_i$'s amongst the $p_i$'s.
\end{rem}

The following corollary gives more precise information on the
equivalence $(i) \Leftrightarrow (iv)$ of the above theorem.

\vspace{5mm}
\begin{cor}
\label{more precise} For an element $f$ in a {\df} $R$, we have $f
\in \mathcal R$ and $\ell(f) = n$ if and only if there exist
$p_1,\dots,p_n \penalty-1000 \in \mathcal A$ such that $R/Rf \cong
\oplus_{i=1}^{n} R/Rp_i$.
\end{cor}
\begin{proof}
If $f\in \mathcal R$ and $\ell(f)=n$ then Lemma \ref{F-basis for
v(f)} implies that there exists an irredundant representation $Rf=
\cap_{i=1}^{n}Rp_i$ with $p_i \in \mathcal A$ and the proof of the
implication $(i) \Rightarrow (iv)$ of the above theorem shows that
$R/Rf \cong \oplus_{i=1}^{n}R/Rp_i$.

  Conversely if $R/Rf \stackrel{\phi}{\longrightarrow}
\bigoplus_{i=1}^nR/Rp_i$ is an isomorphism then using the same
notations as in the proof of the above theorem we have $Rf =
\cap_{i=1}^{n} Rp_i^{x_i}$ and $f\in \mathcal R$.  We must only
show that $n=\ell(f)$.  From Lemma \ref{F-basis for v(f)} this is
equivalent to showing that this representation is irredundant.
Assume at the contrary that this is not the case, without loss of
generality we may assume that $\cap_{i=1}^{n-1} Rp_i^{x_i}
\subseteq Rp_n^{x_n}$.  Now, since $\phi$ is an isomorphism there
exists $h\in R$ such that $\phi(h + Rf) = (0,\dots,0,1 + Rp_n)$
i.e. $hx_i \in Rp_i$ for $i=1,\dots,n-1$ and $hx_n = 1 + Rp_n$.
Since $x_i \notin Rp_i$ we must have $h \in \cap_{i=1}^{n-1}
Rp_i^{x_i}\subseteq Rp_n^{x_n}$.  This implies $hx_n \in Rp_n$ a
contradiction.

\end{proof}
In the following theorem we present different characterizations
for a product to be fully reducible.   Let us first introduce two
relevant definitions :

\vspace{5mm}
\begin{defns} For $a\in R$,
    \begin{enumerate}
      \item[a)] $V'(a):=\{p\in {\mathcal {A}} | \, a\in pR\}$.
      \item[b)] ${\mathbb I}_{R}(Ra)=\{ f \in R | \, af\in Ra \}$.
    \end{enumerate}
\end{defns}

\vspace{5mm}
\begin{thm}
\label{product theorem} For $a,b\in R \setminus U(R)$ the
following are equivalent :
\begin{enumerate}
\item[(i)] $ab$ is fully reducible.
\item[(ii)] $a,b$ are fully reducible and
$R/Rab \cong R/Ra \oplus R/Rb$.
\item[(iii)] $a,b$ are fully reducible and $1\in Ra+bR$.
\item[(iv)] $a,b$ are fully reducible and for all $p\in V(a)$,
$pb$ is fully reducible.
\item[(v)] $a,b$ are fully reducible and for any $F$-basis
$\{p_1,\ldots,p_\ell\}$ of
 $V(a)$, $p_ib$ is fully reducible for $i=1,\ldots,\ell$.
\item[(vi)] $a,b$ are fully reducible and for any $p\in V(a)$ and
any $q\in V'(b)$, $pq$ is fully reducible.
\item[(vii)] $a,b$ are fully reducible and
${\mathbb I}_{R}(Ra)\subset Ra+bR$.
\end{enumerate}
\end{thm}
\begin{proof}
  $(i)\Longrightarrow (ii)$ Since $ab$ is fully
reducible theorem \ref{characterizations of fully reducible} (v)
shows that $a$ and $b$ are also fully reducible. Since $V(b)
\subseteq V(ab)$, we can present an $F$-basis for $V(ab)$ in the
form $p_1,\dots,p_r,q_1,\dots,q_s$ where the $p_i$'s form an
$F$-basis for $V(b)$.  Using Theorem \ref{characterizations of
fully reducible} we obtain the isomorphism $R/Rab \cong
(\oplus_{j=1}^rR/Rp_j) \bigoplus (\oplus_{i=1}^sR/Rq_i) $.  Notice
that we have $r+s = rkV(ab) = \ell (ab) =\ell(a) + \ell (b)$ and
$r=rkV(b)=\ell(b)$ so that $s=\ell(a)=rkV(a)$.  Now, by Lemma
\ref{F-basis for v(f)}, we get $Rb=\cap_{j=1}^{r}Rp_j$ and
$Rab=(\cap_{i=1}^{s}Rq_i) \bigcap (\cap_{j=1}^{r}Rp_j) =
(\cap_{i=1}^{s}Rq_i)\cap Rb = \cap(Rq_i \cap Rb)=\cap Rq_i^{b}b$.
This gives $Ra= \cap _{i=1}^{s}Rq_i^{b}$.   Now, since $\ell(a)=s$
we have (Cf. Remark \ref{remark on thm.}), that $R/Ra \cong \oplus
R/Rq_i^{b} \cong \oplus R/Rq_i$ and $R/Rab \cong (\oplus R/Rq_i)
\bigoplus (\oplus R/Rp_j)\cong R/Ra \oplus R/Rb$.

$(ii)\Longrightarrow (i)$  This is an immediate consequence of
\ref{characterizations of fully reducible} $(iv)$.

  $(i)\Longrightarrow (iii)$ Using the same notations as
in the proof of $(i)\Longrightarrow (ii)$ above,  let us fix an
$F$-basis $\{p_1,\dots,p_r,q_1,\dots,q_s \}$ for $V(ab)$ such that
the $p_i$'s form an $F$-basis for $V(b)$.  We thus have $Rb =
\cap_{j=1}^rRp_j$ and we define $b'$ via $Rb' = \cap_{i=1}^sRq_i$.
Theorem \ref{delta and gamma} then shows that $V(b)\cap V(b') =
\emptyset$ and from lemma \ref{closure,ect} c)
 we get that there exist $u,v\in R$ such that
$ub'+vb=1$.  Left multiplying by $b$, we get $bub'+bvb=b$, in
particular $(bv-1)b\in Rb'$. Therefore $(bv-1)b\in Rb'\cap
Rb=Rab$, and $(bv-1)\in Ra$. This shows that $1\in Ra+bR$ as
desired.

$(iii) \Longrightarrow (iv)$ Let $p$ be an atom in $ V(a)$. Since
$1\in Ra+bR$, $1\in Rp+bR$.  Then there exist $ u,v\in R$ such
that $up+bv=1$. Notice that this shows that $bv\not\in Rp$ and
hence $b\notin Rp^{v}$. Left multiplying $up+bv=1$ by $p$, we get
$pup+pbv=p$. So $pbv\in Rp \cap Rbv$ and hence $pbv \in Rp^{bv}bv$
and $p\in Rp^{bv}$. Since $p$ is an atom we conclude that we must
have $Rp=Rp^{bv}=R(p^{v})^{b}$ where the last equality comes from
\ref{a^b} (b).  We finally get $Rpb = R(p^{v})^{b}b = Rp^{v} \cap
Rb$ which shows that $pb$ is fully reducible since $b$ is fully
reducible.
\newline $(iv) \Longrightarrow (v)$ is obvious.
\newline $(v) \Longrightarrow (i)$ let $A$ be an $F$-basis for
$V(a)$, say $A=\{a_1,\ldots,a_n\}$.  By hypothesis, $\exists
b_1,\ldots,b_n \in {\mathcal {A}}$ such that $Rb_i\cap Rb=Ra_i b$
; let $B=\{b_1,\ldots,b_n\}$. We have $Rab=(\cap Ra_{i})b=\cap
Ra_{i}b= \cap(Rb_i\cap Rb) =\cap Rb_i\cap Rb$. Since $b$ is fully
reducible, this shows that $ab \in \mathcal R$.

  $(iv)\Longrightarrow (vi)$ is obvious and $(vi)\Longrightarrow (iv)$
follows from \ref{characterizations of fully reducible} $(vii)$.

  $(iii)\Longrightarrow (vii)$ Assume $1\in Ra + bR$ and let
  $c\in {\mathbb I}_{R}(Ra)$. We have $ac\in Ra$,
and so $$c=1c\in (Ra+bR)c\subseteq Rac+bR\subseteq Ra+bR$$ as
desired.

  $(vii) \Longrightarrow (iii)$ is trivial since $1\in
{\mathbb I}_{R}(Ra)$.
\end{proof}

\vspace{10mm}

\section{Rank Theorems}
%=======================================================================

\vspace{7mm}

In this final short section we will give some formulas for
computing the rank of algebraic sets of atoms.  Let us first
recall from \ref{full} that a subset $\De \subseteq \mathcal A$ is
full (in $\mathcal A$) if any atom which is $F$-dependent on $\De$
is already in $\De$.

\vspace{5mm}
\begin{prop}
\label{rank of intersection}
 Let $\{\Delta_j : j\in J \}$ be full algebraic sets of atoms.  Then
\begin{enumerate}
\item[a)]$\cap_{j\in J}\Delta_j$ is full algebraic.
\item[b)] If there exists $j_0 \in J$ such that $\Delta_{j_0}$ is
of finite rank then $\cap_{j \in J}\Delta_j$ is of finite rank and

$$ R(\cap_{j\in J}\Delta_j)_\ell=\sum_{j\in J}R(\Delta_j \cap
\Delta_{j_0} )_\ell.$$

%a left least common multiple is given by a
%right greatest common divisor of $\{ (\Delta_j )_\ell , j\in J \}$.
\end{enumerate}
\end{prop}

\begin{proof}
a) Since $\emptyset = \overline{\emptyset}$ is a full algebraic
set, we may assume that $\Delta:=\cap_{j\in J}\Delta_j \ne
\emptyset$.  If $p \in \mathcal A$ is $F$-dependent on $\Delta$
then $p$ is $F$-dependent on each $\Delta_j$ and hence $p \in
\Delta_j$. So $p\in \bigcap_{j\in J} \Delta_j = \Delta$.  This
shows that $\Delta$ is a full algebraic set.

b) Obviously, for any $j\in J, \; \Delta_j \cap \Delta_{j_0}$ is
algebraic of finite rank and since $\cap_j \Delta_j =
\cap_j(\Delta_j \cap \Delta_{j_0})$, we may assume that in fact
all the $\Delta_j$'s are algebraic of finite rank and full ( by a)
above).  Let us put $\Delta = \cap_{j\in J}\Delta_j$.   Let $f\in
R$ and for $j\in J$, let $f_j \in R $ be such that $Rf =
\cap_{\delta \in \Delta} R\delta$ and $Rf_j = \cap_{\delta \in
\Delta_j} R\delta$.  We must show that $R\Delta_\ell=\sum_{j\in
J}R(\Delta_j)_\ell$, i.e. $Rf = \sum_{j \in J}Rf_j$.  Since, for
$j \in J$, $\Delta \subseteq \Delta_j$, we have $Rf_j \subseteq
Rf$.  On the other hand if $h \in R$ is such that $\sum_{j \in
J}Rf_j =Rh$ we have

$$V(h)\subset \bigcap_{j\in J}V(f_j)=\bigcap_{j\in J}\Delta_j =
\Delta = V(f)$$ since $\Delta_j$'s and $\Delta$ are full algebraic
sets.  Therefore, Theorem \ref{characterizations of fully
reducible} again implies that $h$ is a right divisor of $f$.  This
shows that $Rf \subseteq Rh$ and we conclude $\sum_{j \in J}Rf_j
=Rf$, as desired.
\end{proof}

The next theorem gives more precise information than Theorem
\ref{delta and gamma}.

\vspace{5mm}
\begin{thm}
\label{rank of union in general} For any algebraic set of atoms
$\Delta$ and $\Gamma$, we have $$rk(\Delta)+rk(\Gamma)=rk(\Delta
\cup \Gamma)+rk(\overline{\Delta} \cap \overline{\Gamma}).$$
\end{thm}
\begin{proof}
If $\Delta$ or $\Gamma$ is of infinite rank the formula is clear.
We may thus assume that both $\Delta$ and $\Gamma$ are of finite
rank.  Let us write $Rp=R\Gamma_\ell +R\Delta_\ell$ and
$Rq=R\Gamma_\ell \cap R\Delta_\ell$.  Then $Rq=R(\Delta \cup
\Gamma)_\ell$ by \ref{delta and gamma}$(i)$ and
$Rp=R(\overline{\Delta} \cap \overline{\Gamma})_\ell $ by the last
proposition.  Theorem \ref{length of gcd and llcm} gives then
$$\ell(\Gamma_\ell)+\ell(\Delta_\ell)=\ell(q) +\ell(p).$$  In other
words, $$rk(\Delta)+rk(\Gamma)=rk(\Delta \cup
\Gamma)+rk(\overline{\Delta} \cap \overline{\Gamma}).$$
\end{proof}

In order to express the rank of $V(ab)$, let us introduce the
following set : for $a\in R$ we define $I_a := \{q\in {\mathcal
{A}}\ | \ \exists p \in \mathcal {A} \ ; 0\not=  Rp\cap Ra =Rqa
\}$.  Let us also recall our notations : $Rp \cap Ra = Rp^aa$.

\vspace{5mm}
\begin{thm}
Let $a,b \in R$.  Then $$rkV(ba)=rkV(a)+rk(I_{a}\cap V(b)) .$$ In
particular $rkV(ba)\leq rkV(b)+rkV(a)$.
\end{thm}
\begin{proof}
If $V(a) = V(ba)$ we claim that $I_a \cap V(b) = \emptyset$.
Indeed assume $q \in I_a \cap V(b)$, then there exists $p \in
\mathcal A$ and $ a' \in R$ such that $0 \not = Ra \cap Rp = Rqa =
Ra'p$.  In particular there exists $u \in U(R)$ such that $qa
=ua'p$.  Since $q \in V(b)$ we can write $b = b'q$ for some $b'
\in R$.  Multiplying by a on the right gives $ba=b'qa=b'ua'p$.
This shows that $p \in V(ba)=V(a)$ and  hence $Ra \subseteq  Rp$.
We thus get $0 \not = Rqa = Ra \cap Rp = Ra$ and finally $q \in
U(R)$, this is the required contradiction. We may thus assume that
the inclusion $V(a)\subset V(ba)$ is proper.  Let
$\{a_1,\ldots,a_n \}$ be an $F$-basis for $V(a)$ and extend it
into an $F$-basis for $V(ba)$, say
$\{a_1,\ldots,a_n,b_1,\ldots,b_m \}$.  For $i\in \{1,\dots,m\}$ we
 have that $a\notin Rb_i$ and $ba\in Rb_i$.  Then by
 Lemma \ref{product formula},
$b\in Rb_{i}^a$ and by Proposition \ref{bases and independence},
$\{b^{a}_1,\ldots,b^{a}_m \}$ is $F$-independent.  This shows that
$rkV(ba)\leq rkV(a)+rk(I_{a}\cap V(b))$.  For the other inequality
let $\{b^{a}_1,\ldots,b^{a}_m \}$ be an $F$-basis for $I_{a}\cap
V(b)$.  Then by Proposition \ref{bases and independence}
$\{a_1,\ldots,a_n,b_1,\ldots,b_m \}\subset V(ba)$ is
$F$-independent.  This shows that $rkV(ba)\geq rkV(a)+rk(I_{a}\cap
V(b))$.
%Assume $rkV(b)<m$. By the
%lemma above, we know that $\{b^{h}_1,\ldots,b^{h}_m \}$ is  {\it
%F}-independent. But, as $Rba\subset Rb_j$ and $Ra\not\subset
%Rb_j$, by \ref{a^b}, $Rb\subset Rb^{a}_j$ for $j=1,\ldots,m$. So
%$rkV(b)<m\leq rkV(b)$. This gives a contradiction and so this
%shows that $m<rkV(b)$.
% Hence $rkV(ba)\leq m+n\leq rkV(b)+rkV(a)$ as desired.
\end{proof}

\vspace{10mm}

{\centerline {\bf ACKNOWLEDGEMENT}}

\vspace{3mm}

We would like to thank the referee and T.Y.Lam for many helpful
remarks and suggestions.   Thanks to them we avoided awkward flaws
and missprints.

%\section{localizations}

%the definition of a $2-fir$ suggests that there could be some
%relations between these rings and Ore localizations. In this
%section we will investigate these relations and study also some
%other kind of localizations for $2-firs$.

%***what happens if we consider infinite intersection in this thm.
%? Question:  assume f = gh is fully reducible can we always say
%that there exists g' similar to g which divides f on the right,
%and similarly that there always exists an h' similar to h that
%divides f on the left ?***
%$\mathbb{N}$.

\end{document}